\documentclass[12pt,reqno]{amsart}

\usepackage{amscd}
\usepackage{amssymb}
\usepackage{epsfig}
\usepackage{float}
\usepackage{url}
\usepackage{MnSymbol}
\newcommand{\Z}{{\mathbb Z}}
\newcommand{\Q}{{\mathbb Q}}
\newcommand{\C}{{\mathbb C}}
\newcommand{\R}{{\mathbb R}}
\renewcommand{\P}{{\mathbb P}}

\newcommand{\OO}{{\mathcal O}}

\newcommand{\www}{\widetilde}

\newcommand{\oooo}{\overline}
\newcommand{\uuuu}{\underline}

\newcommand{\nnn}{\nabla}

\newcommand{\paa}{\partial}

\DeclareMathOperator{\Gr}{Gr}
\DeclareMathOperator{\Hom}{Hom}
\DeclareMathOperator{\id}{id}
\DeclareMathOperator{\Imm}{Im}

\DeclareMathOperator{\lcm}{lcm}
\DeclareMathOperator{\modd}{mod}

\DeclareMathOperator{\Rad}{Rad}

\DeclareMathOperator{\Seif}{Seif}
\DeclareMathOperator{\Sp}{Sp}
\DeclareMathOperator{\Spp}{Spp}

\DeclareMathOperator{\Tr}{Tr}

\begin{document}

\theoremstyle{plain}
\newtheorem{lemma}{Lemma}[section]
\newtheorem{definition/lemma}[lemma]{Definition/Lemma}
\newtheorem{theorem}[lemma]{Theorem}
\newtheorem{proposition}[lemma]{Proposition}
\newtheorem{corollary}[lemma]{Corollary}
\newtheorem{conjecture}[lemma]{Conjecture}
\newtheorem{conjectures}[lemma]{Conjectures}

\theoremstyle{definition}
\newtheorem{definition}[lemma]{Definition}
\newtheorem{withouttitle}[lemma]{}
\newtheorem{remark}[lemma]{Remark}
\newtheorem{remarks}[lemma]{Remarks}
\newtheorem{example}[lemma]{Example}
\newtheorem{examples}[lemma]{Examples}
\newtheorem{notations}[lemma]{Notations}

\title[Real Seifert forms and Steenbrink PMHS] 
{Real Seifert forms and polarizing forms of
Steenbrink mixed Hodge structures}

\author{Sven Balnojan and Claus Hertling}

\address{Sven Balnojan\\
Lehrstuhl f\"ur Mathematik VI, Universit\"at Mannheim, 
Seminargeb\"aude A 5, 6, 68131 Mannheim, Germany}

\email{sbalnoja@mail.uni-mannheim.de}

\address{Claus Hertling\\
Lehrstuhl f\"ur Mathematik VI, Universit\"at Mannheim, 
Seminargeb\"aude A 5, 6, 68131 Mannheim, Germany}

\email{hertling@math.uni-mannheim.de}



\date{December 01, 2017}

\thanks{This work was supported by the DFG grant He2287/4-1
(SISYPH)}

\maketitle

\begin{abstract}
An isolated hypersurface singularity comes equipped with
many different pairings on different spaces,
the intersection form and the Seifert form on the 
Milnor lattice, a polarizing form for a mixed
Hodge structure on a dual space, and a flat pairing on
the cohomology bundle. 
This paper describes them and their relations systematically 
in an abstract setting. We expect applications also in other
areas than singularity theory. A good part of the paper
is elementary, but not well known: 
the classification of irreducible Seifert form pairs,
the polarizing form on the generalized eigenspace with
eigenvalue 1, 
an automorphism from a Fourier-Laplace transformation
which involves the Gamma function and which 
relates Seifert form and polarizing form
and a flat pairing on the cohomology bundle.
New is a correction of a Thom-Sebastiani formula
for Steenbrink's Hodge filtration in the case of singularities.
It uses the Fourier-Laplace transformation.
A special case is a square root of a Tate twist
for Steenbrink mixed Hodge structures.
\end{abstract}

\maketitle

\tableofcontents

\setcounter{section}{0}

\section{Introduction}\label{s1}
\setcounter{equation}{0}

\noindent
One subject of this paper are real Seifert forms.
In section \ref{s2} a real Seifert form is simply a 
nondegenerate bilinear form $L:H_\R\times H_\R\to\R$
on a finite dimensional $\R$-vector space $H_\R$. The form $L$ is 
in general neither symmetric nor antisymmetric.
It induces an automorphism $M:H_\R\to H_\R$, which is 
called its monodromy, by 
\begin{eqnarray}
L(Ma,b)=L(b,a)\qquad\textup{for }a,b\in H_\R,
\end{eqnarray}
a symmetric bilinear form $I_s:=L+L^t$ (where $L^t(a,b):=L(b,a)$),
and an antisymmetric bilinear form $I_a:=L^t-L$. 
Because of $L(Ma,Mb)=L(Mb,a)=L(a,b)$, the forms $L,I_s$ and 
$I_a$ are $M$-invariant.

The paper starts with this basic linear algebra setting and
studies and enhances it in four steps. The first three steps
are of an abstract nature, and we expect them to have
many different applications, especially in algebraic geometry.
The fourth step is an application to isolated hypersurface
singularities, which is our personal motivation for
developing the material in this paper.

Step 1 = section \ref{s2} shows that any Seifert form pair
$(H_\R,L)$ splits 
(in general not uniquely, but uniquely up to isomorphism)
into an orthogonal and direct sum of irreducible Seifert form
pairs (theorem \ref{t2.5} (b)), it classifies the irreducible
Seifert form pairs (theorem \ref{t2.9}) and gives the 
signatures of their symmetric forms $I_s$ (lemma \ref{t2.10}).

This uses the relation with the notion of an isometric triple
$(H_\R,M,S)$: $H_\R$ is as above, $M$ is an automorphism of
$H_\R$, and $S$ is a nondegenerate $M$-invariant and symmetric
or antisymmetric bilinear form on $H_\R$. 
Milnor \cite[\S 3]{Mi69} classified such triples over 
arbitrary fields. Also an isometric triple splits into an 
orthogonal and direct sum of irreducible pieces 
(theorem \ref{t2.5} (a)). Theorem \ref{t2.8} specializes
Milnor's general classification results to a classification
of irreducible (real) isometric triples.
The lemmata \ref{t2.3} and \ref{t2.4} allow to move between
Seifert form pairs and isometric triples, although the relation
is not one-to-one.

Step 2 = the sections \ref{s3} and \ref{s4} connects 
the Seifert form pairs and the isometric triples with
Steenbrink polarized mixed Hodge structures
(Steenbrink PMHS), an enhancement
of mixed Hodge structures which we define in definition
\ref{t3.3}. Section \ref{s3} reviews mixed Hodge structures,
several enhancements by automorphisms and/or polarizing forms,
and Steenbrink's notions of spectral pairs and spectral numbers
(definition \ref{t3.6}) of Steenbrink mixed Hodge structures.
Theorem \ref{t3.8} gives the irreducible isometric triples
in a Steenbrink PMHS.

Section \ref{s4} connects this with Seifert forms.
It defines a Seifert form $L^{nor}$ 
for a Steenbrink PMHS (definition 
\ref{t4.2} (c)) and gives the irreducible Seifert form
pairs in a Steenbrink PMHS (theorem \ref{t4.4}).
This theorem recovers also a result of Nemethi \cite{Ne95},
namely that the spectral pairs modulo $2\Z\times\{0\}$
are equivalent to the Seifert form of a Steenbrink PMHS.
The sections \ref{s2} to \ref{s4} have some overlap with
the paper \cite{Ne95}. Though he does not consider the
full polarization of a Steenbrink MHS, and he 
classifies explicitly hermitian Seifert form pairs, but
not real Seifert form pairs. In section \ref{s2} we found
it easier to derive the classification or irreducible
Seifert form pairs directly from \cite[\S 3]{Mi69} than
via \cite{Ne95}.

A new ingredient which is neither in \cite{Ne95} nor
in any other papers except \cite{He03}, is an automorphism
(definition \ref{t4.2} (a))
\begin{eqnarray}\label{1.2}
G&=&\bigoplus_{\alpha\in (0,1]}G^{(\alpha)}:H_\C\to H_\C
\quad\textup{with}\\
G^{(\alpha)}&:=&\Gamma(\alpha\cdot\id-\frac{N}{2\pi i})
:H_{e^{-2\pi i\alpha}}\to H_{e^{-2\pi i\alpha}}.\nonumber
\end{eqnarray}
Its definition requires only a finite dimensional
complex vector space $H_\C$ with an automorphism $M$
with eigenvalues $\lambda\in S^1$, semisimple part $M_s$,
unipotent part $M_u$, nilpotent part $N=\log M_u$ and generalized
eigenspaces $H_\lambda:=\ker (M_s-\lambda\cdot\id)$. 
Here $\Gamma(.)$ is the Gamma function.
The true meaning of $G$ 
becomes transparent only in section \ref{s5}
where it arises in a Fourier-Laplace transformation.
But already theorem \ref{t4.3} gives formulas which connect
the polarizing form $S$ and the Seifert form $L^{nor}$ of a 
Steenbrink PMHS with the help of $G$.
We believe that this automorphism $G$ deserves more 
attention than it has obtained up to now.

Step 3 = section \ref{s5} works with a holomorphic vector
bundle on $\C^*$ with a flat holomorphic connection. 
It recalls the well known definition
of elementary sections, the spaces $C^\alpha$ which they
form, and the Malgrange-Kashiwara $V$-filtration.
Not so well known, but elementary is a correspondence
in lemma \ref{t5.1} between three data:
sums of two isometric triples, Seifert form pairs,
and holomorphic bundles on $\C^*$ with a flat holomorphic
connection and a flat real subbundle
and a certain flat pairing $P$ between the fibers at $z\in\C^*$
and $-z$. Theorem \ref{t5.2} enhances this correspondence
with formulas which express a Fourier-Laplace transformation
between elementary sections using $G$ and which 
connect the pairings $P$ and $L^{nor}$. Theorem \ref{t5.2}
and theorem \ref{t4.3} give a relation between $P$ and $S$,
which was stated without proof in \cite[Proposition 7.7]{He03}.

Step 4 = section \ref{s6} is our application to the case
of singularities. An isolated hypersurface singularity
(short: singularity) is a holomorphic function germ
$f:(\C^{m+1},0)\to (\C,0)$ with an isolated singularity at 0.
It had been studied by Milnor \cite{Mi68} and then by
a growing community of singularity theory people.
We recall basic topological notions around it, 
its Milnor lattice $Ml(f)\cong\Z^\mu$, where $\mu\in\Z_{\geq 1}$
is its Milnor number, and on the Milnor lattice its Seifert
form $L$, its monodromy $M$, and its intersection form.
A standard reference is \cite{AGV88}.
An important holomorphic datum of a singularity $f$ is its
Brieskorn lattice $H_0''(f)$, the germ at 0 of a canonical
extension to 0 of its flat cohomology bundle on a punctured
disk $\Delta^*$ \cite{Br70}. 
Varchenko observed that $H_0''(f)$  gives rise to
a MHS \cite{Va80}. Scherk and Steenbrink \cite{SS85}
and M. Saito \cite{SaM82} modified this observation to a recipe
to obtain Steenbrink's MHS $F^\bullet_{St}$ 
from $H_0''(f)$ and the
$V$-filtration. In \cite{He99}\cite{He02} this was enhanced
with a polarizing form $S$  to a (signed) Steenbrink PMHS.
The Fourier-Laplace transformation $FL(H_0''(f))$ and the
relation with the pairing $P$ were considered in \cite{He03}.
Here we recall these facts. The theorems \ref{t4.3}
and \ref{t5.2} are relevant. The tuple 
$TEZP(f):=(Ml(f),m,L,P,H_0''(f))$ can be called a TEZP-structure
(cf. \cite[definition 2.12]{He03}). 

The new point in section \ref{s6} is a Thom-Sebastiani
formula for the TEZP-structures of singularities.
Thom-Sebastiani formulas connect data of a singularity
$f(x_0,...,x_m)$ and data of a singularity 
$g(x_{m+1},...,x_{m+n+1})$ with the data of the singularity $f+g$.
Thom-Sebastiani formulas for Milnor lattice, monodromy and
Seifert form are classical.
Theorem \ref{t6.4} gives the Thom-Sebastiani formula
$TEZP(f)\otimes TEZP(g)\cong TEZP(f+g)$.
An application is a correction of a Thom-Sebastiani formula
in \cite{SS85} for the Hodge filtration $F^\bullet_{St}$
of Steenbrink's MHS. One has to replace in that formula
$F^\bullet_{St}$ by $G(F^\bullet_{St})$. So, here again
the automorphism $G$ is important.

The special case of a suspension, i.e. $f$ as above 
and $g$ with $g=x_{m+1}^2$, leads to a formula which can be
seen as a square root of a Tate twist for Steenbrink PMHS.
It is already stated in theorem \ref{t4.6}. It uses $G$.

This paper collects many classicial pieces. Especially, 
parts are close to \cite{Ne95} and to \cite{SS85}. 
But the relations between the many different pairings, 
the classifications of Seifert forms in general,
their appearance in Steenbrink PMHS,
and the relevance of the automorphism $G$ 
have not been made so explicit before.
We expect applications also in other contexts than
singularities, namely in Landau-Ginzburg models and
in derived algebraic geometry.

\section{Isometric structures and real Seifert forms}\label{s2}
\setcounter{equation}{0}

\noindent
Here Seifert form pairs and isometric triples are defined
and studied. General results in \cite[\S 3]{Mi69} are used
for the classification of isometric triples. This and
their relationship to Seifert form pairs is used for
the classification of Seifert form pairs.
Nemethi \cite{Ne95} undertook the classification of
hermitian Seifert form pairs. One can derive the classification
of real Seifert form pairs from his paper. But we found it
easier to use \cite[\S 3]{Mi69} directly.

\begin{notations}\label{t2.1}
In this section, $H_K$ is a finite dimensional vector space
over a field $K$. If $H_\R$ is given, then $H_\C=H_\R\otimes_\R\C
=H_\R\oplus iH_\R$ is the complexification of $H_\R$.

If $L:H_K\times H_K\to K$ is a bilinear form then
two subspaces $V_1,V_2\subset H_K$ are {\it $L$-orthogonal}
if $L(V_1,V_2)=L(V_2,V_1)=0$.

If $M:H_K\to H_K$ is an automorphism, then
$M_s,M_u,N:H_K\to H_K$ denote its semisimple,
its unipotent and its nilpotent part with
$M=M_sM_u=M_uM_s$ and $N=\log M_u,e^N=M_u$.
If $K=\C$, denote 
$H_\lambda:=\ker(M_s-\lambda\cdot \id):H_\C\to H_\C$,
$H_{\neq 1}:=\bigoplus_{\lambda\neq 1}H_\lambda$,
$H_{\neq -1}:=\bigoplus_{\lambda\neq -1}H_\lambda$.
\end{notations}

\begin{definition}\label{t2.2}
(a) A {\it Seifert form pair} is a pair $(H_\R,L)$ where
$L:H_\R\times H_\R\to\R$ is a nondegenerate bilinear form.
It is called {\it irreducible} if $H_\R$ does not split
into two nontrivial (i.e. both $\neq\{0\}$) 
$L$-orthogonal subspaces.

(b) An {\it isometric triple} is a triple $(H_\R,M,S)$ where
$M:H_\R\to H_\R$ is an automorphism called {\it monodromy},
$S:H_\R\times H_\R\to\R$ is a nondegenerate and 
(for some $m\in\{0,1\}$) $(-1)^m$-symmetric bilinear form
and $M$ is an isometry of $S$.
The triple is called {\it irreducible} if $H_\R$ does not split
into two nontrivial  $S$-orthogonal and $M$-invariant subspaces.
\end{definition}

The following two lemmata show that one can go from Seifert 
form pairs to isometric triples and vice versa,
though the relation is not 1-1.
Starting with $(H_\R,L)$, one has a fixed monodromy $M$ on $H_\R$,
but there are several possible choices of a suitable subspace
$H_\R'$ and a bilinear form $S$ such that $(H_\R',M,S)$ is
an isometric triple. 
Below $I_s$ and $I_a$ are most prominent, but 
$I_s^{(2)},I_a^{(2)},I_s^{(3)}$ and $I_a^{(3)}$ play a role
in the PMHS's of isolated hypersurface singularities.

\begin{lemma}\label{t2.3}
A Seifert form pair $(H_\R,L)$ 
comes equipped with the following data.

(a) Its {\rm monodromy} $M:H_\R\to H_\R$ is the unique 
automorphism with
\begin{eqnarray}\label{2.1} 
L(Ma,b)=L(b,a) \qquad\textup{for all }a,b\in H_\R.
\end{eqnarray}

(b) Define bilinear forms $I_s$ and $I_a$ on $H_\R$,
$I_s^{(2)}$ on $H_\R\cap H_{\neq -1}$, 
$I_a^{(2)}$ on $H_\R\cap H_{\neq 1}$,
$I_s^{(3)}$ on $H_\R\cap H_1$ and $I_a^{(3)}$ on $H_\R\cap H_{-1}$
by
\begin{eqnarray}\label{2.2}
I_s(a,b)&:=& L(b,a)+L(a,b)=L((M+\id)a,b),   \\
I_a(a,b)&:=& L(b,a)-L(a,b)=L((M-\id)a,b),   \nonumber\\
I_s^{(2)}(a,b)&:=& L(a,\frac{1}{M+\id}b)
=I_s(\frac{1}{M+\id}a,\frac{1}{M+\id}b),   \nonumber\\
I_a^{(2)}(a,b)&:=& L(a,\frac{1}{M-\id}b)
=I_a(\frac{1}{M-\id}a,\frac{1}{M-\id}b),  \nonumber\\
I_s^{(3)}(a,b)&:=& L(a,\frac{N}{M-\id}b),   \nonumber\\
I_a^{(3)}(a,b)&:=& L(a,\frac{N}{M+\id}b), \nonumber
\end{eqnarray}
where $\frac{N}{M-\varepsilon\id}$ on 
$H_\R\cap H_\varepsilon$ for $\varepsilon\in\{\pm 1\}$
is the inverse of the automorphism
\begin{eqnarray}\label{2.3}
\frac{M-\varepsilon\id}{N}
&:=& \frac{\varepsilon e^N-\varepsilon\id}{N}:=
\varepsilon\cdot \sum_{k=1}^{\dim H_\R}\frac{1}{k!}\cdot N^{k-1}.
\end{eqnarray}
(Remark that for example in the case $N=0$
$\frac{M-\varepsilon\id}{N}=\varepsilon\id.$)

The bilinear forms $I_s,I_s^{(2)}$ and $I_s^{(3)}$ are
symmetric, the bilinear forms 
$I_a,I_a^{(2)}$ and $I_a^{(3)}$ are antisymmetric.
$I_s^{(2)},I_a^{(2)},I_s^{(3)}$ and $I_a^{(3)}$ are
nondegenerate (on their respective definition domains).
The radical of $I_s$ is $\ker(M+\id)\subset H_{-1}$,
so $I_s$ is nondegenerate on $H_{\neq -1}$.
The radical of $I_a$ is $\ker(M-\id)\subset H_{1}$,
so $I_a$ is nondegenerate on $H_{\neq 1}$.

The automorphisms $M,M_s$ and $M_u$ are isometries of
$L,I_s,I_a,I_s^{(2)},I_a^{(2)},I_s^{(3)}$ and $I_a^{(3)}$,
and $N$ is an infinitesimal isometry of them.
\end{lemma}

{\bf Proof:} 
(a) $M$ is well defined and unique because $L$ is nondegenerate.

(b) $M$ is an isometry of $L$ because applying two times
\eqref{2.1} gives
\begin{eqnarray*}
L(Ma,Mb)=L(Mb,a)=L(a,b).
\end{eqnarray*}
$I_s^{(3)}$ is symmetric and $I_a^{(3)}$ is antisymmetric
because for $\varepsilon\in\{\pm 1\}$ and $a,b\in H_\varepsilon$
\begin{eqnarray*}
L(a,\frac{N}{M-\varepsilon\id}b)
&=&L(M\frac{N}{M-\varepsilon\id}b,a)
=\varepsilon L(\frac{-N}{M^{-1}-\varepsilon\id}b,a)\\
&=&\varepsilon L(b,\frac{N}{M-\varepsilon\id}a).
\end{eqnarray*}
The rest is elementary linear algebra.\hfill$\Box$

\begin{lemma}\label{t2.4}
From an isometric triple one can obtain in different ways
a Seifert form pair. Let $\delta\in\{\pm 1\}$.
Let $(H_\R,M,S)$ be an isometric triple with $S$
$\delta$-symmetric and $H_{-\delta}=\{0\}$,
so $H=H_{\neq -\delta}$ and 
$M+\delta\id$ is invertible. 
Define the Seifert forms $L^{(1)}$ and $L^{(2)}$ by
\begin{eqnarray}\label{2.4}
L^{(1)}(a,b)&:=& S(\frac{1}{M+\delta\id}a,b),\\
L^{(2)}(a,b)&:=& S(a,(M+\delta\id)b).\nonumber
\end{eqnarray}
If $H=H_\delta$, define the Seifert form $L^{(3)}$ by
\begin{eqnarray*}
L^{(3)}(a,b)&:=& S(a,\frac{M-\delta\id}{N}b).
\end{eqnarray*}
For any of these Seifert forms, the monodromy $M$ in lemma
\ref{t2.3} (a) is the monodromy $M$ here. The following table says
which bilinear form in lemma \ref{t2.3} (b) is the $S$ here.
\begin{eqnarray}\label{2.5}
\begin{array}{l|l|l|l|l}
 & L^{(1)} & L^{(2)} & L^{(3)} &  \\ \hline
\delta=1 & I_s & I_s^{(2)} & I_s^{(3)} & =S\\
\delta=-1 & I_a & I_a^{(2)} & I_a^{(3)} & =S\\
\end{array}
\end{eqnarray}
\end{lemma}

{\bf Proof:}
$M$ here and $M$ in lemma \ref{t2.3} (a) coincide because
the $M$ here satisfies
\begin{eqnarray*}
L^{(1)}(Ma,b)&=& S(\frac{M}{M+\delta\id}a,b)
=\delta\cdot S(b,\frac{M}{M+\delta\id}a)\\
&=& \delta\cdot S(\frac{M^{-1}}{M^{-1}+\delta\id}b,a)
=S(\frac{\id}{\delta\id+M}b,a)=L^{(1)}(b,a),
\end{eqnarray*}
and similarly $L^{(2)}(Ma,b)=L^{(2)}(b,a)$,
$L^{(3)}(Ma,b)=L^{(3)}(b,a)$.
The table follows from comparison of the formulas in lemma 
\ref{t2.3} (b) and in lemma \ref{t2.4}.
\hfill$\Box$

\bigskip
Because in a Seifert form pair $(H_\R,L)$
and in an isometric triple $(H_\R,M,S)$, 
the monodromy $M$ is an isometry,
the subspace $H_\lambda$ is $L$-dual respectively
$S$-dual to $H_{\lambda^{-1}}$ and $L$-orthogonal
respectively $S$-orthogonal to all subspaces $H_\kappa$
with $\kappa\neq\lambda^{-1}$. Therefore $H_\R$ splits
canonically into the $M$-invariant and $L$-orthogonal
respectively $S$-orthogonal summands
\begin{eqnarray}\label{2.6}
&&H_\R\cap H_1,\quad H_\R\cap H_{-1},\\
&&H_\R\cap (H_\lambda\oplus H_{\oooo\lambda})\qquad
\textup{for }\lambda\in\{\zeta\in S^1\, |\, \Imm\zeta>0\},
\label{2.7}\\
&&H_\R\cap (H_\lambda\oplus H_{\lambda^{-1}})\qquad\textup{for }
\lambda\in\R_{>1}\cup\R_{<-1},\label{2.8}\\
&&H_\R\cap (H_\lambda\oplus H_{\lambda^{-1}}\oplus 
H_{\oooo\lambda}\oplus H_{\oooo\lambda^{-1}})\label{2.9}\\
&&\hspace*{4cm}\textup{for }\lambda\in\{\zeta\in\C\, |\, |\zeta|>1,
\Imm\zeta>0\}.\nonumber
\end{eqnarray}
In the case of a Seifert form pair, one can choose on
each of these summands a bilinear form $S$ in lemma 
\ref{t2.3} (b) such that $(\textup{the summand},M,S)$
becomes an isometric triple.
Then a splitting of this summand into (irreducible)
Seifert form pairs is a splitting into (irreducible)
isometric triples and vice versa.

Milnor classified isometric triples in \cite[\S 3]{Mi69}
and proved part (a) of the following theorem.
Part (b) is a consequence of part (a) and the lemmata
\ref{t2.3} and \ref{t2.4}.

\begin{theorem}\label{t2.5}
(a) Any isometric triple splits into a direct sum (the summands
are $S$-orthogonal and $M$-invariant) of irreducible
isometric triples. The splitting is unique up to isomorphism.

(b) Any Seifert form pair splits into a direct sum (the summands
are $L$-orthogonal) of irreducible Seifert form pairs. 
The splitting is unique up to isomorphism.
\end{theorem}

It rests to classify the irreducible isometric triples and
via this the irreducible Seifert form pairs.
The irreducible isometric triples had been classified
in an implicit way in \cite[\S 3]{Mi69}.
Nemethi \cite{Ne95} classified the {\it hermitian
Seifert form pairs}, building on \cite[\S 3]{Mi69},
and one can derive from \cite{Ne95} also the irreducible
real Seifert form pairs. But we will use \cite[\S 3]{Mi69}
directly.
We start by examples which in fact will contain all
irreducible isometric triples.

\begin{examples}\label{t2.6}
(i) For $n\in\Z_{\geq 1}$, the following $n\times n$-matrices
will be useful.
\begin{eqnarray*}
E_n=
\begin{pmatrix}1& & \\ & \ddots & \\ & & 1\end{pmatrix},\ 
J_n=(\delta_{j,k+1})_{j,k=1,...,n}=
\begin{pmatrix} 0&&& \\ 1&\ddots&& \\ &\ddots&\ddots& \\&&1&0 
\end{pmatrix},\\
E_n^{per}=((-1)^{j-1}\delta_{j,n+1-k})_{j,k=1,...,n}=
\begin{pmatrix}&&&1\\&&-1&\\&\udots&&\\(-1)^{n-1}&&&\end{pmatrix}.
\end{eqnarray*}

(ii) Choose $n\in\Z_{\geq 1}$, $\lambda\in\{\pm 1\}$ 
and $\varepsilon\in\{\pm 1\}$.
Let $\dim H_\R=n$, and let $\uuuu{a}=(a_1,...,a_n)$ be a 
basis of $H_\R$. Then the monodromy $M$ and the
$(-1)^{n-1}$-symmetric pairing $S$ with
\begin{eqnarray}\label{2.10}
M_s=\lambda\cdot\id,\quad N\uuuu{a}=\uuuu{a}\cdot J_n,\quad
S(\uuuu{a}^t,\uuuu{a})=\varepsilon\cdot E_n^{per}
\end{eqnarray}
give an isometric triple $(H_\R,M,S)$, which is called
$\Tr(\lambda,1,n,\varepsilon)$.
It is irreducible because the monodromy has only one 
Jordan block.

\medskip
(iii) Choose $n\in\Z_{\geq 1}$, $\lambda\in S^1$, 
$\varepsilon\in\{\pm 1\}$ and $m\in\{0,1\}$.
Let $\dim H_\R=2n$.
Choose a complex subspace $H^{(1)}\subset H_\C$
such that $H_\C=H^{(1)}\oplus \oooo{H^{(1)}}$.
Let $\uuuu{a}=(a_1,...,a_n)$ be a basis of $H^{(1)}$.
Then $\uuuu{\oooo{a}}=(\oooo{a}_1,...,\oooo{a}_n)$ is a
basis of $\oooo{H^{(1)}}$. Then the monodromy $M$ and the 
$(-1)^{m}$-symmetric pairing $S$ with
\begin{eqnarray}\label{2.11}
M_s=\lambda\cdot\id|_{H^{(1)}}\oplus
\oooo{\lambda}\cdot\id|_{\oooo{H^{(1)}}},
\quad N\uuuu{a}=\uuuu{a}\cdot J_n,
\quad N\uuuu{\oooo a}=\uuuu{\oooo a}\cdot J_n,\\
S(\begin{pmatrix}\uuuu{a}^t\\ \uuuu{\oooo a}^t\end{pmatrix},
(\uuuu{a},\uuuu{\oooo a}))= i^{n+m+1}\cdot\varepsilon\cdot 
\begin{pmatrix}0 & E_n^{per}\\ (-1)^{n+m+1}E_n^{per} &
0\end{pmatrix} \nonumber
\end{eqnarray}
give an isometric triple $(H_\R,M,S)$, which is called
$\Tr(\lambda,2,n,m,\varepsilon)$. 
Using the basis $(\uuuu{\oooo a},\uuuu{a})$ instead of the
basis $(\uuuu{a},\uuuu{\oooo a})$, one finds
\begin{eqnarray}\label{2.12}
\Tr(\lambda,2,n,m,\varepsilon)\cong
\Tr(\oooo\lambda,2,n,m,(-1)^{n+m+1}\varepsilon).
\end{eqnarray}
If $\lambda\neq\pm 1$ it is irreducible because the
two generalized eigenspaces $H^{(1)}$ and $\oooo{H^{(1)}}$
are $S$-dual (that they are complex conjugate, serves 
equally well) and the monodromy has on each of them only 
one Jordan block.
For $\lambda=\pm 1$ see lemma \ref{t2.7}.

\medskip
(iv) Choose $n\in\Z_{\geq 1}$, $\lambda\in\R_{>1}\cup 
\R_{<-1}$ and $m\in\{0,1\}$.
Let $\dim H_\R=2n$.
Choose a splitting $H_\R= H^{(1)}\oplus H^{(2)}$ into
two $n$-dimensional subspaces.
Let $\uuuu{a}=(a_1,...,a_n)$ be a basis of $H^{(1)}$,
and let $\uuuu{b}=(b_1,...,b_n)$ be a basis of $H^{(2)}$.
Then the monodromy $M$ and the 
$(-1)^{m}$-symmetric pairing $S$ with
\begin{eqnarray}\label{2.13}
M_s=\lambda\cdot\id|_{H^{(1)}}\oplus
\lambda^{-1}\cdot\id|_{H^{(2)}},
\quad N\uuuu{a}=\uuuu{a}\cdot J_n,
\quad N\uuuu{b}=\uuuu{b}\cdot J_n,\\
S(\begin{pmatrix}\uuuu{a}^t\\ \uuuu{b}^t\end{pmatrix},
(\uuuu{a},\uuuu{b}))= 
\begin{pmatrix}0 & E_n^{per}\\ (-1)^{n+m+1}E_n^{per} &
0\end{pmatrix} \nonumber
\end{eqnarray}
give an isometric triple $(H_\R,M,S)$, which is called
$\Tr(\lambda,2,n,m)$. It is irreducible because the
two generalized eigenspaces $H^{(1)}$ and $H^{(2)}$
are $S$-dual and the monodromy has on each of them only 
one Jordan block.

\medskip
(v) Choose $n\in\Z_{\geq 1}$, $\lambda\in\{\zeta\in \C\, |\, 
|\zeta|>1,\Imm\zeta> 0\}$, $\varepsilon\in\{\pm 1\}$ and
$m\in\{0,1\}$. 
Let $\dim H_\R=4n$.
Choose two $n$-dimensional 
complex subspaces $H^{(1)},H^{(2)}\subset H_\C$
such that $H_\C=H^{(1)}\oplus H^{(2)}\oplus 
\oooo{H^{(1)}}\oplus \oooo{H^{(2)}}$.
Let $\uuuu{a}=(a_1,...,a_n)$ be a basis of $H^{(1)}$,
and let $\uuuu{b}=(b_1,...,b_n)$ be a basis of $H^{(2)}$. 
Then the monodromy $M$ and the 
$(-1)^{m}$-symmetric pairing $S$ with
\begin{eqnarray}\label{2.14}
M_s=\lambda\cdot\id|_{H^{(1)}}\oplus
\lambda^{-1}\cdot\id|_{H^{(2)}}\oplus
\oooo{\lambda}\cdot\id|_{\oooo{H^{(1)}}},
\oooo{\lambda}^{-1}\cdot\id|_{\oooo{H^{(2)}}},\hspace*{1cm}\\
N\uuuu{a}=\uuuu{a}\cdot J_n,
\quad N\uuuu{b}=\uuuu{b}\cdot J_n,
\quad N\uuuu{\oooo a}=\uuuu{\oooo a}\cdot J_n, 
\quad N\uuuu{\oooo b}=\uuuu{\oooo b}\cdot J_n, 
\nonumber\hspace*{1cm}\\
S(\begin{pmatrix}\uuuu{a}^t\\ \uuuu{b}^t\\ \uuuu{\oooo a}^t \\
\uuuu{\oooo b}^t\end{pmatrix},
(\uuuu{a},\uuuu{b},\uuuu{\oooo a},\uuuu{\oooo b})) = 
\begin{pmatrix}0 & E_n^{per} & 0 & 0  \\ 
(-1)^{n+m+1}E_n^{per} & 0 &  0 & 0 \\
0 & 0 & 0 & E_n^{per} \\
0 & 0 & (-1)^{n+m+1}E_n^{per} & 0 \end{pmatrix} \nonumber
\end{eqnarray}
give an isometric triple $(H_\R,M,S)$, which is called
$\Tr(\lambda,4,n,m)$. 
It is irreducible because the monodromy has on each of the
four generalized eigenspaces only one Jordan block,
$H^{(2)}$ is $S$-dual to $H^{(1)}$,
$\oooo{H^{(1)}}$ is the complex conjugate of $H^{(1)}$
and $\oooo{H^{(2)}}$ is $S$-dual to $\oooo{H^{(1)}}$
(and the complex conjugate of $H^{(2)}$).
\end{examples}

\begin{lemma}\label{t2.7}
Consider $\lambda\in\{\pm 1\}$.
The types $\Tr(\lambda,1,n,\varepsilon)$ in the examples
\ref{t2.6} (ii) are irreducible and pairwise non-isomorphic.
If $n+m+1\equiv 1(2)$ then by \eqref{2.12}
\begin{eqnarray}\label{2.15}
\Tr(\lambda,2,n,m,1)&\cong& \Tr(\lambda,2,n,m,-1).
\end{eqnarray}
This type is irreducible. If $n+m+1\equiv 0(2)$ then 
$\Tr(\lambda,2,n,m,1)$ and $\Tr(\lambda,2,n,m,-1)$
are not isomorphic and are reducible, 
\begin{eqnarray}\label{2.16}
Tr(\lambda,2,n,m,\varepsilon)
&\cong& 2\cdot \Tr(\lambda,1,n,(-1)^{\frac{n+m+1}{2}}\varepsilon).
\end{eqnarray}
\end{lemma}

{\bf Proof:}
The $\varepsilon$ in $\Tr(\lambda,1,n,\varepsilon)$
is an invariant of the isomorphism class because
$S(b,N^{n-1}b)\in\varepsilon\cdot \R_{>0}$
for any $b\in H_\R-\Imm N$.
Therefore the $\Tr(\lambda,1,n,\varepsilon)$ are
pairwise non-isomorphic.

Now we turn to the examples \eqref{2.6} (iii).
For the proof of \eqref{2.16}, 
work with the real basis $(\uuuu{a}+\uuuu{\oooo a},
i(\uuuu{a}-\uuuu{\oooo a}))$.
One has to calculate the matrix of $S$
for the new basis. Details are left to the reader.

Irreducibility of $\Tr(\lambda,2,n,m,1)$ in the case
$n+m+1\equiv 1(2)$:
Indirect proof. Suppose $H_\R=V_1\oplus V_2$ is an
$S$-orthogonal and $M$-invariant splitting.
Then each of $V_1$ and $V_2$ consists of one Jordan block of $M$.
Choose a basis $\uuuu{c}=(c_1,...,c_n)$ of $V_1$
with $N\uuuu{c}=\uuuu{c}\cdot J_n$. 
Use that $S$ is here $(-1)^n$-symmetric and that $N$
is an infinitesimal isometry. It gives
\begin{eqnarray*}
S(c_j,c_{n+1-j})&\stackrel{N}{=}&
(-1)^{(n+1-j)-j}\cdot S(c_{n+1-j},c_j) \stackrel{S}{=}
-S(c_j,c_{n+1-j}),\textup{ so }=0,\\
S(c_j,c_{n+1-k})&\stackrel{N}{=}& 0\quad \textup{for }k<j\textup{ anyway.}
\end{eqnarray*}
Then $S$ is degenerate on $V_1$, a contradiction.
\hfill $\Box$

\begin{theorem}\label{t2.8}
\cite[\S 3]{Mi69}
The irreducible isometric triples are given by the following
types, which are all non-isomorphic.
\begin{eqnarray}\label{2.17}
\Tr(\lambda,1,n,\varepsilon)&\textup{with}&
\lambda\in\{\pm 1\},\\
\Tr(\lambda,2,n,m,1)&\textup{with}&
\lambda\in\{\pm 1\}\ \&\ m\equiv n(2),\label{2.18}\\
\Tr(\lambda,2,n,m,\varepsilon)&\textup{with}&
\lambda\in\{\zeta\in S^1\, |\, \Imm\zeta>0\},\label{2.19}\\
\Tr(\lambda,2,n,m)&\textup{with}&
\lambda\in\R_{>1}\cup \R_{<-1},\label{2.20}\\
\Tr(\lambda,4,n,m)&\textup{with}&
\lambda\in\{\zeta\in \C\, |\, |\zeta|>1,\Imm\zeta>0\}.\label{2.21}
\end{eqnarray}
Here $n\in\Z_{\geq 1},\varepsilon\in\{\pm 1\},
m\in\{0,1\}$.
\end{theorem}

{\bf Proof:}
As this is only implicit in \cite[\S 3]{Mi69},
we provide additional arguments.

The cases $\lambda\in\R_{>1}\cup\R_{<-1}$
and $\lambda\in \{\zeta\in \C\, |\, |\zeta|>1,\Imm\zeta>0\}$
are subsumed in \cite[\S 3]{Mi69} as  \lq\lq case 3\rq\rq  \ 
and are the easiest cases.
Consider $\lambda\in \{\zeta\in \C\, |\, |\zeta|>1,\Imm\zeta>0\}$,
and consider an isometric triple
$(H_\R,M,S)$ with 
$H_\C=H_\lambda\oplus H_{\lambda^{-1}}\oplus H_{\oooo\lambda}
\oplus H_{\oooo{\lambda}^{-1}}$ and 
$S$ $(-1)^m$-symmetric. Choose a basis
$\uuuu{a}=(a_1,...,a_n)$ of $H_\lambda$ which is adapted
to the Jordan block structure of $N$ on $H^{(1)}$, so
\begin{eqnarray*}
\uuuu{a}=(\uuuu{a}^{(1)},...,\uuuu{a}^{(r)})\textup{ with }
N\uuuu{a}^{(j)}=\uuuu{a}^{(j)}\cdot J_{n_j}
\quad\textup{for some }r,n_1,...,n_r\in\Z_{\geq 1}
\end{eqnarray*}
(so $n_1+...+n_r=n$). 
Let $\uuuu{c}=(\uuuu{c}^{(1)},...,\uuuu{c}^{(r)})$
be the $S$-dual basis of $H_{\lambda^{-1}}$.
Define 
\begin{eqnarray*}
\uuuu{b}^{(j)}&:=&((-1)^{n_j-1}b_{n_j}^{(j)},
(-1)^{n_j-2}b_{n_j-1}^{(j)},...,
-b_2^{(j)},b_1^{(j)})\quad\textup{and}\\
\uuuu{b}&:=&(\uuuu{b}^{(1)},...,\uuuu{b}^{(r)}).
\end{eqnarray*}
Then $H_\C$ splits into the $S$-orthogonal and 
$M$-invariant subspaces 
$\langle \uuuu{a}^{(j)},\uuuu{b}^{(j)},
\uuuu{\oooo a}^{(j)},\uuuu{\oooo b}^{(j)}\rangle$
for $j=1,...,r$, and the $j$-th space is with this  basis
of the type $\Tr(\lambda,4,n_j,m)$.

The case $\lambda\in \R_{>1}\cup \R_{<-1}$ is similar.

\medskip
The cases $\lambda\in\{\zeta\in S^1\, |\, \Imm\zeta>0\}$
and $\lambda=\pm 1$ are called 
\lq\lq case 1\rq\rq  \ respectively \lq\lq case 2\rq\rq  
\ in \cite[\S 3]{Mi69}.
For such a value $\lambda$ let $(H_\R,M,S)$ be an isometric
triple with $S$ $(-1)^m$-symmetric for some
$m\in\{0,1\}$ 
and with $H_\C=H_\lambda\oplus H_{\oooo\lambda}$ in
case 1 and $H_\C=H_\lambda$ in case 2.

Theorem 3.2 in \cite{Mi69} says that the isometric triple
splits into isometric triples such that on each summand
all Jordan blocks have the same length
and that the summands are unique up to isomorphism.
Therefore suppose that on $H_\C$ all Jordan blocks have
the same length $n$.

Now consider first case 1, so 
$\lambda\in\{\zeta\in S^1\, |\, \Imm\zeta>0\}$.
The sesquilinear (=linear$\times$semilinear) form $S_{res,1}$ on
$H_\lambda/(H_\lambda\cap\Imm N)$ with
\begin{eqnarray}\label{2.22}
S_{res,1}([a],[b])&:=&
(-i)^{n+m+1}\cdot S(a,N^{n-1}\oooo{b})
\quad\textup{for }a,b\in H_\lambda
\end{eqnarray}
is well defined and nondegenerate and hermitian:
It is well defined and nondegenerate because $N$ is an 
infinitesimal isometry and all Jordan blocks have the
same length $n$, so that especially
$\ker N=\Imm N^{n-1}$ and $S(\Imm N,\ker N)=0$.
The following calculation shows that it is hermitian,
\begin{eqnarray*}
S_{res,1}([b],[a])&=& 
(-i)^{n+m+1}\cdot S(b,N^{n-1}\oooo{a})\\
&=&(-1)^{n+m+1}(-1)^m\cdot S(N^{n-1}\oooo{a},b)\\
&=& (-i)^{n+m+1}(-1)^{n+m+1}\oooo{S(a,N^{n-1}\oooo{b})}\\
&=&\oooo{S_{res,1}([a],[b])}.
\end{eqnarray*}
Theorem 3.3 in \cite{Mi69} implies that the isomorphism class
of the isometric triple $(H_\R,M,S)$ is determined by the
signature of $S_{res,1}$.

In the case $\Tr(\lambda,2,n,m,\varepsilon)$ we have 
$H_\lambda/(H_\lambda\cap \Imm N)=\C\cdot[a_1]$
and 
\begin{eqnarray*}
S_{res,1}([a_1],[a_1])= (-i)^{n+m+1}\cdot S(a_1,N^{n-1}\oooo{a_1})
=\varepsilon.
\end{eqnarray*}
Therefore in the general case above, 
the isometric triple $(H_\R,M,S)$ is isomorphic
to a sum of triples $\Tr(\lambda,2,n,m,\varepsilon_j)$
for $j=1,2,...,\frac{1}{2n}\dim H_\R$ where the
$\varepsilon_j\in\{\pm 1\}$ are determined by the
signature of $S_{res,1}$.

Finally consider case 2, so $\lambda=\pm 1$.
The bilinear form $S_{res,2}$ on
$H_\R/\Imm N$ with
\begin{eqnarray}\label{2.23}
S_{res,2}([a],[b])&:=&S(a,N^{n-1}b)
\quad\textup{for }a,b\in H_\R
\end{eqnarray}
is well defined and nondegenerate and $(-1)^{n+m+1}$-symmetric:
It is well defined and nondegenerate for the same
reasons as $S_{res,1}$. 
The following calculation shows that it is 
$(-1)^{n+m+1}$-symmetric,
\begin{eqnarray*}
S_{res,2}([b],[a])&=& S(b,N^{n-1}a)
=(-1)^m\cdot S(N^{n-1}a,b)\\
&=& (-1)^{n+m+1}S(a,N^{n-1}b)
=(-1)^{n+m+1}S_{res,2}([a],[b]).
\end{eqnarray*}
Theorem 3.4 in \cite{Mi69} implies that the isomorphism class
of the isometric triple $(H_\R,M,S)$ is determined by the
signature of $S_{res,2}$ if $n+m+1\equiv 0(2)$
and that it is independent of any additional data if
$n+m+1\equiv 1(2)$. 

In the cases $\Tr(\lambda,1,n,\varepsilon)$
with $\lambda=\pm 1$ and $n+m+1\equiv 0(2)$ we have 
$H_\R/\Imm N=\R\cdot[a_1]$
and 
\begin{eqnarray*}
S_{res,2}([a_1],[a_1])= S(a_1,N^{n-1}a_1)=\varepsilon.
\end{eqnarray*}
Therefore in the general case above, 
the isometric triple $(H_\R,M,S)$ is 
in the case $n+m+1\equiv 0(2)$ isomorphic
to a sum of triples $\Tr(\lambda,1,n,\varepsilon_j)$
for $j=1,2,...,\frac{1}{n}\dim H_\R$ where the
$\varepsilon_j\in\{\pm 1\}$ are determined by the
signature of $S_{res,2}$.
In the case $n+m+1\equiv 1(2)$, the isometric triple
$(H_\R,M,S)$ is isomorphic to a sum of triples
$\Tr(\lambda,2,n,m,1)$.
\hfill$\Box$ 

\bigskip

Theorem \ref{t2.8} together with the lemmata \ref{t2.3}
and \ref{t2.4} gives also the classification of the
irreducible Seifert form pairs in theorem \ref{t2.9}.
The proof of theorem \ref{t2.9} states which
isometric triples give rise to which Seifert form pairs.

\begin{theorem}\label{t2.9}
The irreducible Seifert form pairs are given by the types
with the following names.
\begin{eqnarray}\label{2.24}
\Seif(\lambda,1,n,\varepsilon)&\textup{with}&
(\lambda=1\ \& \ n\equiv 1(2))\\
&\textup{or}& (\lambda=-1\ \&\ n\equiv 0(2)),\nonumber\\
\Seif(\lambda,2,n)&\textup{with}&
(\lambda=1\ \& \ n\equiv 0(2))\label{2.25}\\
&\textup{or}& (\lambda=-1\ \&\ n\equiv 1(2)),\nonumber\\
\Seif(\lambda,2,n,\zeta)
&\cong& \Seif(\oooo{\lambda},2,n,\oooo{\zeta})
\label{2.26}\\
&\textup{with}& \lambda,\zeta\in S^1-\{\pm 1\},
\zeta^2=\oooo{\lambda}\cdot (-1)^{n+1},\nonumber\\
\Seif(\lambda,2,n)&\textup{with}&
\lambda\in\R_{>1}\cup \R_{<-1},\label{2.27}\\
\Seif(\lambda,4,n)&\textup{with}&
\lambda\in\{\zeta\in \C\, |\, |\zeta|>1,\Imm\zeta>0\}.\label{2.28}
\end{eqnarray}
Here $n\in\Z_{\geq 1},\varepsilon\in\{\pm 1\}$.
The types are uniquely determined by the properties 
above of $\lambda$ and $n$ and the following properties.

\begin{list}{}{}
\item[\eqref{2.24}] $\Seif(\lambda,1,n,\varepsilon):$ 
$\dim H_\R=n$, $H_\C=H_\lambda$, 
one Jordan block, for each $a\in H_\R-\Imm N$
$$L(a,N^{n-1}a)\in \varepsilon\cdot\R_{>0}.$$
\item[\eqref{2.25}] $\Seif(\lambda,2,n):$
$\dim H_\R=2n$, $H_\C=H_\lambda$, 
two Jordan blocks of size $n$. 
\item[\eqref{2.26}] $\Seif(\lambda,2,n,\zeta):$
$\dim H_\R=2n$, $H_\C=H_\lambda\oplus H_{\oooo\lambda}$, 
two Jordan blocks, for each $a\in H_\lambda-\Imm N$
$$L(a,N^{n-1}\oooo{a})\in \zeta\cdot\R_{>0}.$$
\item[\eqref{2.27}] $\Seif(\lambda,2,n):$
$\dim H_\R=2n$, $H_\C=H_\lambda\oplus H_{\lambda^{-1}}$, 
two Jordan blocks of size $n$. 
\item[\eqref{2.28}] $\Seif(\lambda,4,n):$
$\dim H_\R=4n$, $H_\C=H_\lambda\oplus H_{\lambda^{-1}}
\oplus H_{\oooo\lambda}\oplus H_{\oooo{\lambda}^{-1}}$, 
four Jordan blocks of size $n$.
\end{list}
\end{theorem}

{\bf Proof:}
The following table lists irreducible isometric triples
and chosen Seifert forms from lemma \ref{t2.4} which give
rise to irreducible Seifert form pairs.
In the cases \eqref{2.24} and \eqref{2.26}, calculations 
after the table 
show that the Seifert form pairs have the stated properties.
In all cases \eqref{2.24}--\eqref{2.28}, one sees that the
stated properties characterize the Seifert form pairs
uniquely by going back via lemma \ref{t2.3} to 
isometric triples and comparing their classification
in theorem \ref{t2.8}.

Lemma \ref{t2.4} will be applied now.
The $\delta$ in lemma \ref{t2.4} is here
in the table in the case \eqref{2.24} 
$\delta=\lambda=(-1)^{n-1}$, 
in the other cases $\delta=(-1)^m$.
\begin{eqnarray*}
\begin{array}{llll}
 & & L\textup{ from lemma \ref{t2.4}} & \\
\eqref{2.24} & \Tr(\lambda,1,n,\varepsilon) & 
L^{(1)}\textup{ or }L^{(2)}\textup{ or }L^{(3)} &
\Seif(\lambda,1,n,\lambda\cdot\varepsilon) \\
\eqref{2.25} & \Tr(\lambda,2,n,m,1) & 
L^{(1)}\textup{ or }L^{(2)}\textup{ or }L^{(3)} &
\Seif(\lambda,2,n) \\
& \textup{with }m\equiv n(2) & & \\
\eqref{2.26} & \Tr(\lambda,2,n,m,\varepsilon) & 
L^{(1)}\textup{ or }L^{(2)} &
\Seif(\lambda,2,n,
\frac{\oooo{\lambda}+1}{|\lambda+1|}i^{n+1}\varepsilon) \\
\eqref{2.27} & \Tr(\lambda,2,n,m) & 
L^{(1)}\textup{ or }L^{(2)} &
\Seif(\lambda,2,n) \\
\eqref{2.28} & \Tr(\lambda,4,n,m) & 
L^{(1)}\textup{ or }L^{(2)} &
\Seif(\lambda,4,n) 
\end{array}
\end{eqnarray*}
The calculation for the case \eqref{2.24} with 
$L^{(1)}$ ($L^{(2)}$ and $L^{(3)}$ are analogous):
\begin{eqnarray*}
L^{(1)}(a,N^{n-1}a) &=& 
S(\frac{1}{M+\lambda\id}a,N^{n-1}a) \\
&=& S(\frac{1}{2\lambda}a(+\textup{something in }\Imm N),
N^{n-1}a)\\
&=& \frac{1}{2}\lambda\cdot S(a,N^{n-1}a)=\frac{1}{2}
\cdot\lambda\cdot \varepsilon.
\end{eqnarray*}
The calculation for the case \eqref{2.26} with 
$L^{(2)}$ ($L^{(1)}$ is analogous):
\begin{eqnarray*}
L^{(2)}(a,N^{n-1}a) &=& 
S(a,(M+(-1)^m\id)N^{n-1}\oooo{a}) \\
&=& S(a,(\oooo{\lambda}+(-1)^m)N^{n-1}\oooo{a})\\
&=& (\oooo{\lambda}+(-1)^m)\cdot i^{n+m+1} \cdot\varepsilon\\
&\in& (\oooo{\lambda}+1)\cdot i^{n+1}\cdot\varepsilon
\cdot \R_{>0}.
\end{eqnarray*}
In the last line $\Imm(\lambda)>0$ (in \eqref{2.19} for
$\Tr(\lambda,2,n,m,\varepsilon)$) is used.
\hfill$\Box$

\bigskip
The next lemma gives for each irreducible Seifert form pair
the signature of $I_s$.
This is useful if one wants to determine the irreducible
pieces  of a given Seifert form pair.
Here the signature $(p,q,r)$ means
$p:=\max(\dim U\, |\, U\textup{ pos. def. subspace of }H_\R)$,
$q:=\dim\Rad I_s$, 
$r=n-p-q=\max(\dim U\, |\, U\textup{ neg. def. subspace of }H_\R)$.

\begin{lemma}\label{t2.10}
The following table lists for the irreducible Seifert form
pairs in theorem \ref{t2.9} the signature of $I_s$ and for all
cases with $\Rad I_s=\{0\}$ the type of the irreducible
isometric triple. 
\begin{eqnarray*}
\begin{array}{llcl}
\textup{type of a Seifert} & \textup{form pair}
& \textup{signature of }I_s
& \textup{isometric str.}\\ \hline
\Seif(1,1,n,\varepsilon) 
&\textup{with }n\equiv\varepsilon(4)
& (\frac{n+1}{2},0,\frac{n-1}{2}) 
& \Tr(1,1,n,\varepsilon) \\
\Seif(1,1,n,\varepsilon) 
&\textup{with }n\equiv-\varepsilon(4)  
& (\frac{n-1}{2},0,\frac{n+1}{2}) 
& \Tr(1,1,n,\varepsilon) \\
\Seif(-1,1,n,\varepsilon) 
&\textup{with }n-1\equiv\varepsilon(4)
& \hspace*{0.3cm} (\frac{n}{2},1,\frac{n-2}{2})
&  \\
\Seif(-1,1,n,\varepsilon) 
&\textup{with }n-1\equiv-\varepsilon(4)
& (\frac{n-2}{2},1,\frac{n}{2}) \hspace*{0.4cm}
& \\
\Seif(1,2,n) 
&\textup{(with }n\equiv 0(2))
& (n,0,n) 
& \Tr(1,2,n,0,1) \\
\Seif(-1,2,n) &
\textup{(with }n\equiv 1(2))   
& (n-1,2,n-1) 
&  \\
\Seif(\lambda,2,n,\zeta\varepsilon) 
&\textup{with } n\equiv 0(2)
& (n,0,n) 
& \Tr(\lambda,2,n,0,\varepsilon) \\ 
&  (\textup{and }\lambda\in S^1-\{\pm 1\}) & & \\
\Seif(\lambda,2,n,\zeta) 
&\textup{ with }  n\equiv 1(2) 
& (n-1,0,n+1) 
& \Tr(\lambda,2,n,0,1) \\
& (\textup{and }\lambda\in S^1-\{\pm 1\}) & & \\
\Seif(\lambda,2,n,-\zeta) 
&\textup{with }  n\equiv 1(2) 
& (n+1,0,n-1) 
& \Tr(\lambda,2,n,0,-1) \\
& (\textup{and }\lambda\in S^1-\{\pm 1\}) & & \\
\Seif(\lambda,2,n) 
&\textup{with }
\lambda\in \R_{>1}\cup\R_{<-1} 
& (n,0,n) 
& \Tr(\lambda,2,n,0) \\
\Seif(\lambda,4,n) 
&\textup{with }  
\lambda\in\{\zeta\in\C| 
& (2n,0,2n) 
& \Tr(\lambda,4,n,0) \\
& |\zeta|>1,\Imm\zeta>0\} & & 
\end{array}
\end{eqnarray*}
Here $n\in\Z_{\geq 1},\varepsilon\in\{\pm 1\}$, and 
in the lines 7--9 
$\zeta:=\frac{\oooo{\lambda}+1}{|\lambda+1|}\cdot i^{n+1}$.
\end{lemma}

{\bf Proof:}
For all cases except those in the lines 3, 4 and 6,
$(H_\R,M,I_s)$ is an irreducible isometric triple,
and the proof of theorem \ref{t2.9} tells which it is.
Then one can read off the signature of $I_s$ from the
examples \ref{t2.6}.

The least easy cases are in the lines 8 and 9. We treat the 
case in line 9 and leave the other cases to the reader.
The case in line 9 is a special case of example \ref{t2.6} (iii).
Here $I_s$ has the same signature as the hermitian matrix
\begin{eqnarray*}
S(\begin{pmatrix} \uuuu{a}^t\\ \uuuu{\oooo{a}}^t\end{pmatrix},
(\uuuu{\oooo{a}},\uuuu{a})) = 
i^{n-1}\cdot \begin{pmatrix}E^{per}_n & 0 \\ 0 & E^{per}_n
\end{pmatrix}.
\end{eqnarray*}
The signature is $(n+1,0,n-1)$.

In the cases in the lines 3, 4 and 6, lemma \ref{t2.3} says
$\Rad I_s=\ker(M+\id)=\ker N$. The induced isometric triple
$(H_\R/\Rad I_s,M,I_s)$ has eigenvalue $-1$ and in the cases
in the lines 3 and 4 only one Jordan block of size $n-1$
and in the cases in the line 6 two Jordan blocks of sizes
$n-1$. Theorem \ref{t2.5} and \ref{t2.8} tell us:
The isometric triple $(H_\R/\Rad I_s,M,I_s)$ is in all cases
irreducible. It is of the type
$\Tr(-1,1,n-1,\www{\varepsilon})$ with a suitable
$\www{\varepsilon}$ in the lines 3 and 4
and of the type $\Tr(-1,2,n-1,0,1)\cong\Tr(-1,2,n-1,0,-1)$
(with $n-1\equiv 0(2)$) in line 6.
The type $\Tr(-1,2,n-1,0,\pm 1)$ has signature 
$(n-1,0,n-1)$. This gives $(n-1,2,n-1)$ in line 6.

The cases in the lines 3 and 4: 
$\www{\varepsilon}$ has to be determined.
For each $a\in H_\R-\Imm N$ we have 
$L(a,N^{n-1}a)\in\varepsilon\cdot\R_{>0}$. 
\begin{eqnarray*}
&&I_s(a,N^{n-2}a) = L(a,N^{n-2}a)+L(N^{n-2}a,a)\\
&=&2L(a,N^{n-2}a) = 2L(N^{n-2}a,a)
= 2 L(Ma,N^{n-2}a)\\
&=&2L(-e^Na,N^{n-2}a) =-2L(a+Na,N^{n-2}a)\\
&=& -2L(a,N^{n-2}a)+2L(a,N^{n-1}a),\quad\textup{thus it is}\\
&=& L(a,N^{n-1}a)\in \varepsilon\cdot\R_{>0},
\end{eqnarray*}
so $\www\varepsilon=\varepsilon$. 
The signature of $I_s$ on $H_\R/\Rad(I_s)$ is the
signature of $\varepsilon\cdot E^{per}_{n-1}$.
\hfill $\Box$

\bigskip

We finish this section with some elementary statements
on induced structures on the dual space.

\begin{notations}\label{t2.11}
Let $H_K$ be a finite dimensional $K$-vector space.
$H_K^\vee:=\Hom(H_K,K)$ is the dual space, and 
$<,>:H_K^\vee\times H_K\to K$ denotes the natural pairing.
If $M:H_K\to H_K$ is an automorphism, then
$M^\vee:H_K^\vee\to H_K^\vee$ is defined by
$<M^\vee a,Mb>=<a,b>$.
If $L:H_K\times H_K\to K$ is a nondegenerate pairing,
let $L^{lin}:H^\vee_K\to H_K$ be the induced isomorphism
with $L(a,b)=<(L^{lin})^{-1}(a),b>$,
and define $L^\vee:H_K^\vee\times H_K^\vee\to K$ by
$L^\vee(a,b)=<a,L^{lin}b>= L(L^{lin}a,L^{lin}b)$.
\end{notations}

\begin{lemma}\label{t2.12}
(a) If $(H_\R,L)$ is a Seifert form pair with monodromy $m$,
then $(L^{lin})^{-1}:(H_\R,L,M)\to (H_\R^\vee,L^\vee,M^\vee)$
is an isomorphism of Seifert form pairs with monodromies.

(b) If $(H_\R,M,S)$ is an isometric triple,
then $(S^{lin})^{-1}:(H_\R,M,S)\to (H_\R^\vee,M^\vee,S^\vee)$
is an isomorphism of isometric triples.

(c) Let $(H_\R,L)$ be a Seifert form pair with 
$H=H_{\neq -\delta}$ for some $\delta\in\{\pm 1\}$.
Denote $S:=I_s$ if $\delta=1$ and $S:=I_a$ if $\delta=-1$.
Then 
\begin{eqnarray}\label{2.29}
L^{lin}\circ M^\vee =M\circ L^{lin},\\
S^{lin} =L^{lin}\circ \frac{1}{M^\vee+\delta\id}
=\frac{1}{M+\delta\id}\circ L^{lin},\label{2.30}\\
S^\vee =S^{L^\vee,(2)}\quad\textup{with }
S^{L^\vee,(2)} (a,b):= 
L^\vee(a,\frac{1}{M^\vee+\delta\id}),\label{2.31}
\end{eqnarray}
so $S^{L^\vee,(2)}$ is the pairing $I_s^{(2)}$
respectively $I_a^{(2)}$ in lemma \ref{t2.3} (b),
but for $L^\vee$ instead of $L$.
\end{lemma}

{\bf Proof:} Elementary. \hfill$\Box$

\section{Polarized mixed Hodge structures}\label{s3}
\setcounter{equation}{0}

\noindent
Steenbrink defined mixed Hodge structures for isolated
hypersurface singularities and their spectral pairs.
These mixed Hodge structures are special in several aspects.
They come equipped with an automorphism of the 
vector space which induces the 
weight filtration and which is essential for the spectral pairs.
And they come equipped with a natural polarization.
Though the spectral pairs are defined without using the
polarization.

Usually a $\Z$-lattice or a $\Q$-vector space underly
a mixed Hodge structure. They give a rigidity and richness
which are usually precious. But we do not want this
rigidity here, so we will not consider a $\Z$-lattice
or a $\Q$-vector space here.

\begin{notations}\label{t3.1}
The notations \ref{t2.1} will be used again.
All filtrations in this paper are finite and exhaustive.
An upper index means a decreasing filtration,
a lower index means an increasing filtration.
The Gauss bracket is denoted $\lfloor .\rfloor:\R\to\Z$.
The upper Gauss bracket is denoted $\lceil.\rceil:\R\to\Z$.
The following two functions will allow to treat
several cases simultaneously:
\begin{eqnarray*}
[.]_2:\Z\to\{0,1\}&&\textup{with}\quad n\equiv [n]_2\modd 2,\\ 
\theta:S^1\to\{0,1\}&&\textup{with }\theta(1):=1
\textup{ and }\theta(\lambda):=0\textup{ for } 
\lambda\neq 1.
\end{eqnarray*}
\end{notations}

The following lemma from \cite[Lemma 6.4]{Sch73}
(see also e.g. \cite[Lemma 2.1]{He99}) prepares
definition \ref{t3.3}. 
It is stated in \cite{Sch73} with $H_\Q$ instead of $H_\R$.

\begin{lemma}\label{t3.2}
Let $m\in \Z$, 
$H_\R$ a finite dimensional $\R$-vector space,\\
$S:H_\R\times H_\R\to \R$ a nondegenerate 
$(-1)^m$-symmetric bilinear form,
and $N:H_\R\to H_\R$ a nilpotent endomorphism which is an 
infinitesimal isometry of $S$. 

(a) There exists a unique increasing filtration 
$W_\bullet \subset H_\R$ such that
$N(W_l)\subset W_{l-2}$ and such that 
$N^l:\Gr^W_{m+l}\to \Gr^W_{m-l}$ is an isomorphism.
Sometimes it will be called $W^{(N,m)}_\bullet$.

(b) $S(W_k,W_{l})=0$ if $k+l< 2m$.

(c) A nondegenerate $(-1)^{m+l}$-symmetric bilinear form $S_l$ is 
well defined on $\Gr^W_{m+l}$ for $l\geq 0$ by the requirement:
$S_l(a,b)=S(\tilde{a}, N^l \tilde{b})$ if $\tilde{a},\tilde{b} \in W_{m+l}$
represent $a,b\in \Gr^W_{m+l}$.

(d) The primitive subspace $P_{m+l}$ of $\Gr^W_{m+l}$ 
is defined by
$$P_{m+l}=\ker (N^{l+1}:\Gr^W_{m+l}\to \Gr^W_{m-l-2})$$
if $l\geq 0$ and $P_{m+l}=0$ if $l<0$.  Then
$$\Gr^W_{m+l}=\bigoplus_{i\geq 0} N^i P_{m+l+2i},$$
and this decomposition is orthogonal with respect to $S_l$ if $l\geq 0.$
\end{lemma}

\begin{definition}\label{t3.3}
(a) A {\it mixed Hodge structure} (short: MHS)
is a tuple $(H_\R,H_\C,F^\bullet,W_\bullet)$ with
$F^\bullet\subset H_\C$ a decreasing {\it Hodge filtration}
and $W_\bullet\subset H_\R$ an increasing 
{\it weight filtration} such that $F^\bullet \Gr^W_k$ gives
a pure Hodge structure of weight $k$ on $\Gr^W_k$, i.e. 
\begin{eqnarray}\label{3.1}
\Gr^W_k=F^p\Gr^W_k\oplus\oooo{F^{k+1-p}\Gr^W_k}.
\end{eqnarray}

(b) A {\it Steenbrink MHS} of weight $m\in\Z$ is a MHS
$(H_\R,H_\C,F^\bullet,W_\bullet)$ together with an
automorphism $M$ (called {\it monodromy}) 
of $(H_\R,H_\C,W_\bullet)$ with the following properties: 
Its semisimple part maps $F^p$ to $F^p$, its nilpotent part $N$
maps $F^p$ to $F^{p-1}$, and $N$ determines $W_\bullet$
as follows.
\begin{eqnarray}\label{3.2}
W_\bullet|_{H_{\neq 1}}=W_\bullet^{(N,m)}\textup{ on }H_{\neq 1},
\quad \textup{and }
W_\bullet|_{H_1}=W_\bullet^{(N,m+1)}\textup{ on }H_1.
\end{eqnarray}

(c) \cite{CK82}\cite{He99}
A {\it polarized mixed Hodge structure} (short: PMHS)
of weight $m\in\Z$ is a tuple
$(H_\R,H_\C,F^\bullet,W_\bullet,N,S)$
with $(m,H_\R,H_\C,S,N,W_\bullet)$ as in lemma \ref{t3.2} and
\begin{list}{}{}
\item[(i)] $(H_\R,H_\C,F^\bullet,W_\bullet)$ is a MHS.
\item[(ii)] $N(F^p)\subset F^{p-1}$.
\item[(iii)] $S(F^p,F^{m+1-p})=0.$
\item[(iv)] The pure Hodge structure $F^\bullet P_{m+l}$
of weight $m+l$ on $P_{m+l}$ is polarized by $S_l$, i.e.
\begin{list}{}{}
\item[($\alpha$)] 
$S_l(F^pP_{m+l},F^{m+l+1-p}P_{m+l})=0.$
\item[($\beta$)]
$i^{2p-m-l}\cdot S_l(a,\oooo{a})>0\textup{ for }
a\in F^pP_{m+l}\cap\oooo{F^{m+l-p}P_{m+l}}-\{0\}.$
\end{list}
\end{list}

\medskip
(d) A {\it Steenbrink PMHS} of weight $m\in\Z$ is a
Steenbring MHS together with a nondegenerate pairing
$S$ such that the restriction to $H_{\neq 1}$ is a PMHS
of weight $m$ and the restriction to $H_1$ is a 
PMHS of weight $m+1$
(especially, $S$ is $(-1)^m$-symmetric on $H_{\neq 1}$
and $(-1)^{m+1}$-symmetric on $H_1$).
\end{definition}

\begin{remarks}\label{t3.4}
In \cite{CK82} condition (c)(iii) is omitted.
Condition (c)(iii) implies condition (iv)($\alpha$)
(therefore we could have omitted condition (iv)($\alpha$)).
In the case of an
isolated hypersurface singularity, the polarization
on $H_1$ was not considered by Steenbrink,
only later in \cite{He99}.
\end{remarks}

Deligne defined subspaces $I^{p,q}$ of a MHS which 
split the Hodge filtration and the weight filtration
in a natural way \cite{De71}. They also behave
well with respect to morphisms and a polarizing form
\cite{CK82}\cite{He99}.

\begin{lemma}\label{t3.5}
For a MHS define
\begin{eqnarray*}
I^{p,q}&:=& \left(F^p\cap W_{p+q}\right)\cap 
\left({\overline F}^q \cap W_{p+q}+
\sum_{j>0}{\overline F}^{q-j}\cap W_{p+q-j-1}\right).
\end{eqnarray*}
Then
\begin{eqnarray}\label{3.3}
F^p&=&\bigoplus_{i,q:\, i\geq p}I^{i,q},\\
W_l&=&\bigoplus_{p+q\leq l}I^{p,q},\label{3.4}\\ 
I^{q,p}&\cong& \oooo{I^{p,q}}\modd W_{p+q-2}.\label{3.5}
\end{eqnarray}
If $W=W^{(N,m)}$ for a nilpotent endomorphism
$N:H_\R\to H_\R$ and a weight $m\in\Z$, then define
additionally for $p+q\geq m$
\begin{eqnarray*}
I_0^{p,q}&:=&\ker (N^{p+q-m+1}:I^{p,q}\to 
I^{m-q-1,m-p-1}).
\end{eqnarray*}
Then
\begin{eqnarray}\label{3.6}
N(I^{p,q})&\subset& I^{p-1,q-1},\\
I^{p,q}&=&\bigoplus_{j\geq 0}N^jI_0^{p+j,q+j},\label{3.7}\\
I^{q,p}_0&\cong& \oooo{I^{p,q}_0}\modd W_{p+q-2}.\label{3.8}
\end{eqnarray}
In the case of a PMHS of weight $m$ with polarizing form $S$
\begin{eqnarray}\label{3.9}
S(I^{p,q},I^{r,s})&=&0\quad\textup{for }(r,s)\neq (m-p,m-q),\\
S(N^iI_0^{p,q},N^jI_0^{r,s})&=&0\quad\textup{for }
(r,s,i+j)\neq (q,p,p+q-m).\label{3.10}
\end{eqnarray}
\end{lemma}

Steenbrink's spectral pairs provide a very intuitive
picture which allows to see and understand the discrete data
in a Steenbrink MHS well.

\begin{definition}\label{t3.6}
\cite{St77}
Let $(H_\R,H_\C,F^\bullet,W_\bullet,M)$ be a Steenbrink
MHS of weight $m$ with $n:=\dim H_\R\geq 1$.
The {\it spectral pairs} are $n$ pairs 
$(\alpha,k)\in\R\times \Z$ with multiplicities 
$d(\alpha,k)\in\Z_{\geq 0}$,
\begin{eqnarray}
\Spp &=& \sum_{(\alpha,k)}d(\alpha,k)\cdot (\alpha,k)
\in\Z_{\geq 0}[\R\times\Z],\nonumber\\
d(\alpha,k)&:=& \dim\Gr^{\lfloor m-\alpha\rfloor}_F
\Gr^W_{k+\theta(\lambda)} H_\lambda
\quad\textup{for }e^{-2\pi i\alpha}=\lambda\label{3.11}
\end{eqnarray}
($\theta(\lambda)$ was defined in the notations \ref{t3.1}).
The {\it spectral numbers} are the first entries in the
spectral pairs,
\begin{eqnarray}
\Sp &=& \sum_{\alpha}d(\alpha)\cdot (\alpha)
\in\Z_{\geq 0}[\R],\nonumber\\
d(\alpha)&:=& \sum_kd(\alpha,k)
=\dim\Gr^{\lfloor m-\alpha\rfloor}_F H_\lambda
\quad\textup{for }e^{-2\pi i\alpha}=\lambda.\label{3.12}
\end{eqnarray}
\end{definition}

Now we will discuss the geometry in the spectral pairs.
Lemma \ref{3.5} will be crucial.
Consider some $p,q\in\Z$ and $\lambda\in S^1$ such that
the space $(I^{p,q}_0)_\lambda:=I^{p,q}_0\cap H_\lambda$
is not $\{0\}$. Then $p+q=m+\theta(\lambda)+l$ for some
$l\in\Z_{\geq 0}$, and the spaces in the two sequences
\begin{eqnarray}\label{3.13}
(I^{p,q}_0)_\lambda, N(I^{p,q}_0)_\lambda,...,
N^l (I^{p,q}_0)_\lambda, \\
(I^{q,p}_0)_{\oooo\lambda}, N(I^{q,p}_0)_{\oooo\lambda},...,
N^l (I^{q,p}_0)_{\oooo\lambda},\label{3.14}
\end{eqnarray}
have all the same dimension. 
They give rise to the following 
{\it ordered pair of spectral pair ladders}, where each 
spectral pair has the same multiplicity
$\dim (I^{p,q}_0)_\lambda$:
\begin{eqnarray}\label{3.15}
(\alpha,m+l),(\alpha+1,m+l-2),...,(\alpha+l,m-l),\\
(m-l-1-\alpha,m+l),(m-l-\alpha+1,m+l-2),...,\nonumber \\
(m-1-\alpha,m-l). \label{3.16}
\end{eqnarray}
In one row the first entry is increasing by 1,
the second entry is decreasing by 2.
Here $\alpha\in\R$ is determined by 
$e^{-2\pi i\alpha}=\lambda$ and $p=\lfloor m-\alpha\rfloor
=m-\lceil\alpha\rceil$.
The first spectral pair $(\alpha,m+l)$ in 
the first spectral pair ladder \eqref{3.15}
comes from $(I^{p,q}_0)_\lambda$.
The first spectral pair $(m-l-1-\alpha,m+l)$ in the second
spectral pair ladder \eqref{3.16}
comes from $(I^{q,p}_0)_{\oooo\lambda}$, because
$q+p=m+\theta(\lambda)+l$, 
$e^{-2\pi i (m-l-1-\alpha)}=\oooo{\lambda}$, and
\begin{eqnarray*}
\lfloor m-(m-l-1-\alpha)\rfloor &=& l+1+\lfloor \alpha\rfloor
=l+\theta(\lambda)+\lceil\alpha\rceil\\
&=& l+\theta(\lambda)+m-p=q.
\end{eqnarray*}
The other spectral pairs follow from the first ones
by applying \eqref{3.6} repeatedly.

If $(p,\lambda)=(q,\oooo{\lambda})$ 
(so $\lambda\in\{\pm 1\}$) then
$(I^{p,q}_0)_\lambda=(I^{q,p}_0)_{\oooo\lambda}$ and then
there is only one spectral pair ladder,
i.e. \eqref{3.15} and \eqref{3.16} agree and their
multiplicity is $\dim (I^{p,p}_0)_\lambda$.
Then the spectral pair ladder is its own partner.
By \eqref{3.7} $\Spp$ consists completely of 
spectral pair ladders, namely pairs of them and 
(for $(p,\lambda)=(q,\oooo{\lambda})$) single ones.
Each pair of spectral pair ladders and also the single ones
are invariant under the Kleinian group
$\id,\pi_1,\pi_2,\pi_3:\R\times\Z\to\R\times\Z$ with
\begin{eqnarray}\label{3.17}
\pi_1:(\frac{m-1}{2}+\alpha,m+l)&\mapsto&
(\frac{m-1}{2}-\alpha,m-l),\\
\pi_2:(\frac{m-1-l}{2}+\alpha,m+l)&\mapsto&
(\frac{m-1-l}{2}-\alpha,m+l),\nonumber\\
\pi_3=\pi_1\circ\pi_2=\pi_2\circ\pi_1:
(\alpha,m+l)&\mapsto&
(\alpha+l,m-l).\nonumber
\end{eqnarray}
Obviously, the decomposition of $\Spp$ into ordered pairs of 
spectral pair ladders and single ones with these symmetries
is unique up to changing the order of the ordered pairs.
If a spectral pair ladder starts at $(\alpha,m+l)$,
its {\it length} is $l+1$ and the {\it distance} 
to its partner is $2\alpha+l+1-m$.
The single ones have distance 0.
Thus for the single ones
\begin{eqnarray}
\alpha=\frac{m-l-1}{2}\in\frac{1}{2}\Z
\quad\textup{ and }\nonumber\\
l\equiv 0(2)\textup{ if }\alpha\in\frac{m-1}{2}+\Z,
\quad l\equiv 1(2)\textup{ if }\alpha\in\frac{m}{2}+\Z.
\label{3.18}
\end{eqnarray}

The following picture gives an example. Each dot stands for
a spectral pair or the corresponding space
$N^j(I^{p,q}_0)_\lambda$. Dots of the same shape
indicate spectral pairs in one pair of 
spectral pair ladders.
Only in order not to overload the picture, we restrict
in the picture to $\alpha\in\frac{1}{2}\Z$.
Also $F^\bullet$ and $W_\bullet$ can be read off.

\noindent
\includegraphics[width=1.0\textwidth]{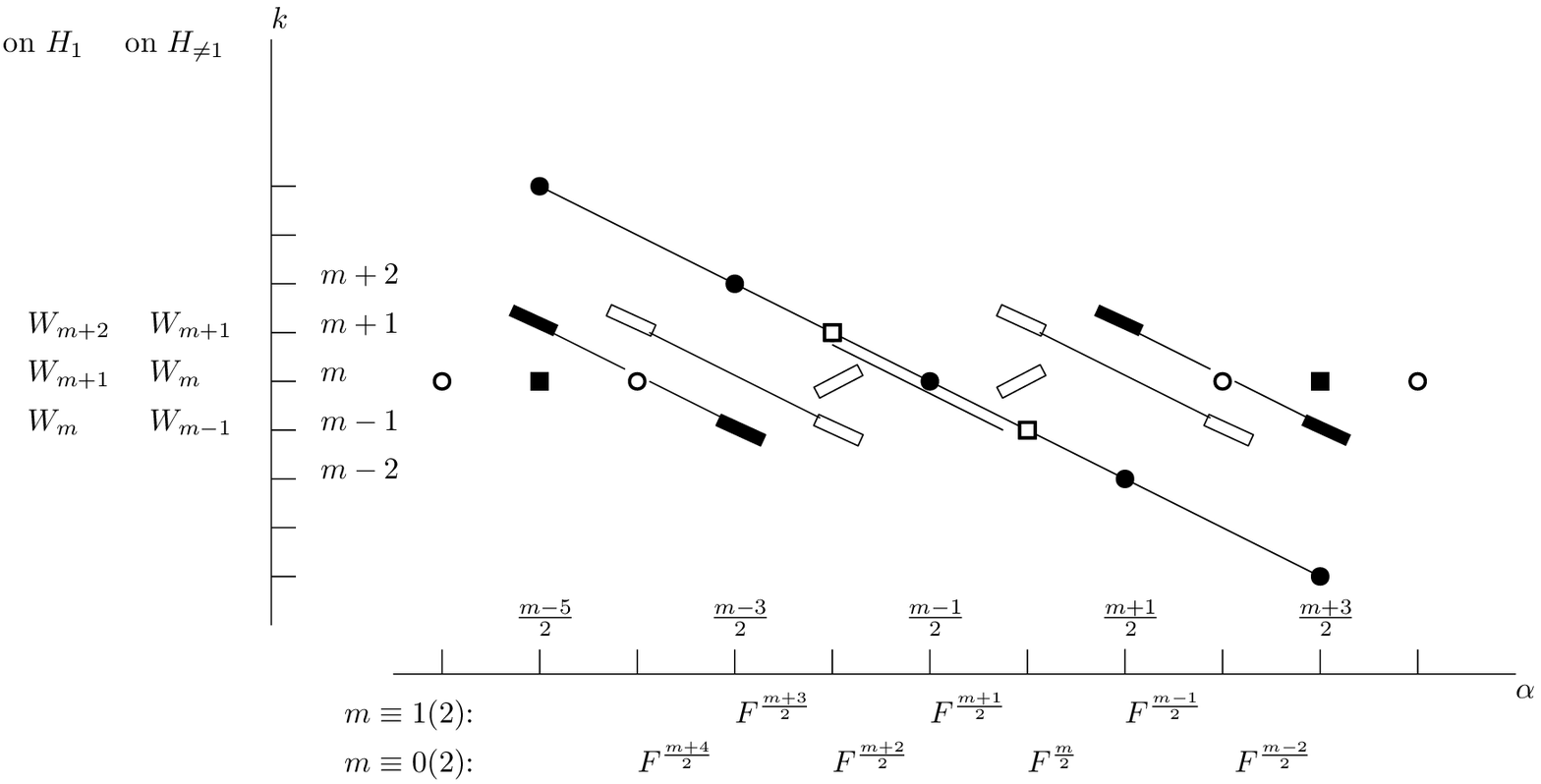} 

In the definition \ref{t3.6} 
of $\Spp$ only a Steenbrink MHS is needed, 
not a Steenbrink PMHS. 
But if we have a Steenbrink PMHS, then it makes sense
to study the underlying isometric structures on
$H_{\neq 1}$ and on $H_1$. 
Theorem \ref{t3.8} studies the isometric triples
$(H_\R\cap H_{\neq 1},M,S)$ and $(H_\R\cap H_1,M,S)$
The following observation from \cite{He99} 
simplifies this study.

\begin{remark}\label{t3.7}
Starting with a reference PMHS $(H_\R,H_\C,F^\bullet_0,
W_\bullet, N,S)$ of some weight $m$,
a classifying space $D_{PMHS}$ for all Hodge filtrations
$F^\bullet$ on $H_\C$ such that $(H_\R,H_\C,F^\bullet,W_\bullet,
N,S)$ is a PMHS with the same spectral pairs as the
reference PMHS was constructed in \cite{He99}.
It contains a filtration $F^\bullet_1$ such that
the $I^{pq}(F^\bullet_1)$ satisfy
$I^{qp}(F^\bullet_1)=\oooo{I^{pq}(F^\bullet_1)}$.
Such a PMHS is called {\it split}.
All this holds also for Steenbrink PMHS.
\end{remark}

\begin{theorem}\label{t3.8}
Let $(H_\R,H_\C,F^\bullet,W_\bullet,M,S)$ be a Steenbrink PMHS
of weight $m$. Because of remark \ref{t3.7} we can suppose
that it is split, i.e. $I^{qp}=\oooo{I^{pq}}$.

(a)  The sum $(H_\R\cap H_{\neq 1},M,S)
\oplus (H_\R\cap H_1,M,S)$
of isometric triples decomposes into the isometric triples
\begin{eqnarray}\label{3.19}
&&\sum_{j=0}^l N^j(I^{p,q}_0)_\lambda + 
\sum_{j=0}^l N^j(I^{q,p}_0)_{\oooo\lambda}\\
&&\textup{for }p,q,\lambda\textup{ with }
(I^{p,q}_0)_\lambda\neq\{0\},\ 
(p,\lambda)\neq(q,\oooo{\lambda}),\ \Imm((-1)^{m+1}\lambda)\geq 0
\nonumber
\end{eqnarray}
and the isometric triples 
\begin{eqnarray}\label{3.20}
\sum_{j=0}^l N^j(I^{p,p}_0)_\lambda\quad 
\textup{for }p,\lambda\textup{ with }\lambda=\pm 1,\ 
(I^{p,p}_0)_\lambda\neq\{0\}. 
\end{eqnarray}

(b) Each of the isometric triples in \eqref{3.19}
decomposes into $\dim (I^{p,q}_0)_\lambda$ many copies
of the isometric triple
\begin{eqnarray}\label{3.21}
\Tr(\lambda,2,l+1,[m+\theta(\lambda)]_2,
(-1)^{\lceil\alpha\rceil-1-\frac{1}{2}(m-\theta(\lambda)+
[m+\theta(\lambda)]_2)}).
\end{eqnarray}
Each of the isometric triples in \eqref{3.20}
decomposes into $\dim (I^{p,p}_0)_\lambda$ many copies
of the isometric triple
\begin{eqnarray}\label{3.22}
\Tr(\lambda,1,l+1,
(-1)^{\lceil\alpha\rceil -\frac{1}{2}(m-\theta(\lambda)-l)}),\\
\textup{and then }(-1)^{m+1}\lambda=(-1)^l.\nonumber
\end{eqnarray}
\end{theorem}

{\bf Proof:}
(a) $I^{qp}=\oooo{I^{pq}}$ implies 
$I^{qp}_0=\oooo{I^{pq}_0}$ and
$(I^{qp}_0)_{\oooo\lambda}=\oooo{(I^{pq}_0)_\lambda}$.
Therefore the spaces in \eqref{3.19} and \eqref{3.20}
are complexifications of real subspaces.

The polarizing form $S$ is $M$-invariant.
The decomposition is $S$-orthogonal by \eqref{3.10}
and the $M$-invariance of $S$. It is obviously
$M$-invariant.

(b) Formula \eqref{3.10} and the $M$-invariance of
$S$ show that the isometric triples
in \eqref{3.19} and \eqref{3.20} are sums of
isometric triples of the types
$\Tr(\lambda,2,l+1,[m+\theta(\lambda)]_2,\varepsilon)$
and $\Tr(\lambda,1,l+1,\varepsilon)$
for suitable $\varepsilon$. 
Here $S$ is $(-1)^{[m+\theta(\lambda)]_2}$-symmetric.
Therefore in the case $\Tr(\lambda,1,l+1,\varepsilon)$
$l\equiv [m+\theta(\lambda)]_2\modd 2$, i.e. 
$(-1)^{m+1}\lambda=(-1)^l$.

It rests to determine
$\varepsilon$. Choose $a\in (I^{p,q}_0)_\lambda-\{0\}$.
The polarizing condition  (c)(iv)($\beta$) in definition
\ref{t3.3} says 
\begin{eqnarray*}
0< i^{p-q}\cdot S_l(a,\oooo{a})
=i^{2p-m-\theta(\lambda)-l}\cdot S(a,N^l\oooo{a}).
\end{eqnarray*}
The following calculations use also $p=m-\lceil\alpha\rceil$.

Consider first the case \eqref{3.19}.
The definition in example \ref{t2.6} (iii) says
\begin{eqnarray*}
S(a,N^l\oooo{a})\in i^{(l+1)+[m+\theta(\lambda)]_2+1}
\cdot\varepsilon\cdot \R_{>0}.
\end{eqnarray*}
Then
\begin{eqnarray*}
\varepsilon &=& i^{-(2p-m-\theta(\lambda)-l)
-(l+1+[m+\theta(\lambda)]_2+1)}\\
&=& (-1)^{-p-1+
\frac{1}{2}(m+\theta(\lambda)-[m+\theta(\lambda)]_2)}\\
&=& (-1)^{\lceil\alpha\rceil -1-
\frac{1}{2}(m-\theta(\lambda)+[m+\theta(\lambda)]_2)}.
\end{eqnarray*}
Consider now the case \eqref{3.20}. The isometric triple
must be one in exampe \ref{t2.6} (ii), so then
$m+\theta(\lambda)\equiv l(2)$ and
$S(a,N^la)\in\varepsilon\cdot \R_{>0}$ and
\begin{eqnarray*}
\varepsilon &=& i^{-(2p-m-\theta(\lambda)-l)}\\
&=& (-1)^{-p+
\frac{1}{2}(m+\theta(\lambda)+l)}
= (-1)^{\lceil\alpha\rceil -
\frac{1}{2}(m-\theta(\lambda)-l)}.
\end{eqnarray*}
\hfill$\Box$

\section{Seifert forms and Steenbrink PMHS}\label{s4}
\setcounter{equation}{0}

\noindent
The purpose of this and the next section is to compare
and relate several bilinear forms: the polarizing
form of a Steenbrink PMHS, a Seifert form and,
in section \ref{s5}, a pairing on a flat bundle
on $\C^*$. They all arise in the case of an 
isolated hypersurface singularity.
But here we consider them abstractly.

Lemma \ref{t4.1} starts with a Seifert form and gives a 
family of together symmetric forms and a hermitian form.

Definition \ref{t4.2} and the theorems \ref{t4.3}
and \ref{t4.4} start from a Steenbrink PMHS.
A (normalized) Seifert form is defined,
and also an automorphism $G$, which seems to have received
less attention than it deserves. Its significance
will become fully transparent only in section \ref{s5}
when a Fourier-Laplace transformation is considered.

Theorem \ref{t4.3} fixes the relations between the 
polarizing form, the Seifert form and this automorphism.
Theorem \ref{t4.4} classifies the irreducible Seifert form
pairs in a Steenbrink PMHS.
It recovers the result of Nemethi \cite{Ne95} that the
spectral pairs $\Spp\modd 2\Z\times\{0\}$ are equivalent
to the Seifert form (and the weight $m$, which we need
as our Seifert form is normalized, but $\Spp$ is not).
Finally, theorem \ref{t4.6} gives for a Steenbrink PMHS
a square root of a Tate twist.
This uses the automorphism $G$.
It is modelled after the suspension of a singularity.

\begin{lemma}\label{t4.1}
Let $(H_\R,L)$ be a Seifert form pair.

(a) For $\lambda$ with $H_\lambda\neq\{0\}$ and $\kappa$
with $\kappa^2=\lambda$ define a pairing
\begin{eqnarray}\label{4.1}
L^{sym}_\kappa:H_\lambda\times H_{1/\lambda}\to\C,\quad
L^{sym}_\kappa(a,b):=\kappa\cdot L(a,e^{-N/2}b).
\end{eqnarray}
Then $L^{sym}_\kappa$ and $L^{sym}_{1/\kappa}$ satisfy
together the symmetry condition
\begin{eqnarray}\label{4.2}
L^{sym}_{1/\kappa}(b,a)=L^{sym}_\kappa(a,b).
\end{eqnarray}

(b) For $\lambda\in S^1$ with $H_\lambda\neq\{0\}$ 
and for $\kappa$ with $\kappa^2=\lambda$ define a sesquilinear
pairing 
\begin{eqnarray}\label{4.3}
L^{herm}_\kappa:H_\lambda\times H_\lambda\to\C,\quad
L^{herm}_\kappa(a,b):=L^{sym}_\kappa(a,\oooo{b}).
\end{eqnarray}
It is hermitian.
\end{lemma}

{\bf Proof:} (a)
\begin{eqnarray*}
L^{sym}_{1/\kappa}(b,a)&=& \kappa^{-1}L(b,e^{-N/2}a)
=\kappa^{-1}L(Me^{-N/2}a,b)\\
&=&\kappa^{-1}L(\lambda e^{N/2}a,b)
=\kappa^{-1}\lambda L(a,e^{-N/2}b)=L^{sym}_\kappa(a,b).
\end{eqnarray*}
(b) Here $\kappa\in S^1$, so $\kappa^{-1}=\oooo\kappa$.
\begin{eqnarray*}
L^{herm}_\kappa(b,a)&=& L^{sym}_\kappa (b,\oooo{a})
=L^{sym}_{1/\kappa}(\oooo{a},b)
=\oooo{\kappa} L(\oooo{a},e^{-N/2}b)\\
&=& \oooo{\kappa L(a,e^{-N/2}\oooo{b})} 
=\oooo{L^{sym}_\kappa(a,\oooo{b})}
=\oooo{L^{herm}_\kappa(a,b)}.
\qquad\Box
\end{eqnarray*}

\begin{definition}\label{t4.2}
Let $(H_\R,H_\C,F^\bullet,W_\bullet,M,S)$ be a Steenbrink
PMHS of weight $m$.

(a) Then each eigenvalue of $M$ is in $S^1$. Define an 
automorphism $G:H_\C\to H_\C$ as follows.
\begin{eqnarray}
G&:=&\bigoplus_{\alpha\in(0,1]}G^{(\alpha)}\textup{ with }
G^{(\alpha)}:H_{e^{-2\pi i\alpha}}\to H_{e^{-2\pi i\alpha}},
\nonumber\\
G^{(\alpha)}&:=& \sum_{k\geq 0}\frac{1}{k!}\Gamma^{(k)}(\alpha)
\left(\frac{-N}{2\pi i}\right)^k=\Gamma(\alpha\cdot \id
-\frac{N}{2\pi i}).\label{4.4}
\end{eqnarray}
$G$ does not respect $H_\R$ if $N\neq 0$. But it commutes
with $M$ and $M_s$ and $N$ and it respects $W_\bullet$.

(b) The {\it normalized} monodromy is $M^{nor}:=(-1)^{m+1}M$,
so $M^{nor}_s=(-1)^{m+1}M_s$, $N^{nor}=N$.
But in (c) and in theorem \ref{t4.3}, $H_\lambda$ and 
$H_{\neq 1}$ still refer to $M$, not to $M^{nor}$.

(c) The {\it normalized} Seifert form 
$L^{nor}:H_\R\times H_\R\to \R$ is defined as follows.
First define an automorphism $\nu:H_\R\to H_\R$ by
\begin{eqnarray}\label{4.5}
\nu&:=&\left\{\begin{array}{ll}
\frac{1}{M-\id}&\textup{on }H_{\neq 1},\\
\frac{-N}{M-\id}&\textup{on }H_1.\end{array}\right. 
\end{eqnarray}
Now define $L^{nor}$ by
\begin{eqnarray}\label{4.6}
L^{nor}(a,b):=
-S(a,\nu^{-1}b)=\left\{\begin{array}{ll}
(-1)^m\cdot L^{(2)}(a,b)&\textup{ on }H_{\neq 1},\\
(-1)^{m+1}\cdot L^{(3)}(a,b)&\textup{ on }H_1.\end{array}\right.
\end{eqnarray}
Here $L^{(2)}$ and $L^{(3)}$ come from lemma \ref{2.4}
and the isometric triples $(H_\R\cap H_{\neq 1},M^{nor},S)$
(with $\delta=(-1)^m$) and $(H_\R\cap H_1,M^{nor},S)$
(with $\delta=(-1)^{m+1}$).
Thus $M^{nor}$ is the monodromy of $L^{nor}$.
\end{definition}

\begin{theorem}\label{t4.3}
Let $(H_\R,H_\C,F^\bullet,W_\bullet,M,S)$ be a Steenbrink
PMHS of weight $m$.

(a) For $a\in H_\lambda,b\in H_{\oooo{\lambda}}$ with
$\lambda=e^{-2\pi i \alpha}$ and $0<\alpha<1$
\begin{eqnarray}\label{4.7}
S(a,b)=\frac{-1}{2\pi i}\cdot e^{-\pi i\alpha}\cdot
L^{nor}(Ga,e^{-N/2}Gb).
\end{eqnarray}
For $a,b\in H_1$
\begin{eqnarray}\label{4.8}
S(a,b)=L^{nor}(Ga,e^{-N/2}Gb).
\end{eqnarray}

(b) $(H_\R,H_\C,G(F^\bullet),W_\bullet,M)$ is a Steenbrink
MHS of weight $m$ with 
\begin{eqnarray}\label{4.9}
L^{sym}_{\kappa\cdot i^{m-1}}(G(F^p)H_\lambda,
G(F^{m+1-p})H_{\oooo{\lambda}})=0&&
\textup{for }\lambda\neq 1\\
&&\textup{ and }\kappa
\textup{ with }\kappa^2=\lambda,\nonumber\\
L^{sym}_{i^{m-1}}(G(F^p)H_1,G(F^{m+2-p})H_1)=0.&&
\label{4.10}
\end{eqnarray}
Remark that $(\kappa\cdot i^{m-1})^2=\lambda(-1)^{m-1}$ 
is the eigenvalue of $M^{nor}$ on $H_\lambda$.

\medskip
(c) Recall the relation between $(I^{pq}_0)_\lambda$
and the first spectral pair $(\alpha,m+l)$ in the
spectral pair ladder in \eqref{3.15}:
$p+q=m+\theta(\lambda)+l$, $p=\lfloor m-\alpha\rfloor$,
$\lambda=e^{-2\pi i \alpha}$. 
Recall also that $m-l-1-\alpha$ is the first spectral number of
the partner spectral pair ladder and that
$2\alpha+l+1-m$ is the distance from the spectral pair ladder
to its partner.

For $a\in (I^{pq}_0)_\lambda-\{0\}$ as well as for 
$a\in G((I^{pq}_0)_\lambda)-\{0\}$
\begin{eqnarray}\label{4.11}
L^{nor}(a,N^l\oooo{a})\in e^{\frac{1}{2}\pi i(2\alpha+l+1-m)}
\cdot\R_{>0}.
\end{eqnarray}
\end{theorem}

{\bf Proof:}
(a) Recall the following identities 
of the Gamma function (they are equivalent if one uses
$\Gamma(x+1)=x\Gamma(x)$).
\begin{eqnarray*}
\Gamma(x)\Gamma(1-x)&=&\frac{\pi}{\sin \pi x}
= e^{\pi ix}\frac{2\pi i}{e^{2\pi i x}-1},\\
\Gamma(1+x)\Gamma(1-x)&=& \frac{\pi x}{\sin \pi x}
=e^{\pi ix}\frac{2\pi i x}{e^{2\pi i x}-1}.
\end{eqnarray*}
They imply for $0<\alpha<1$
\begin{eqnarray*}
\Gamma(\alpha\id+\frac{N}{2\pi i})
\Gamma(\id-(\alpha\id+\frac{N}{2\pi i}))
&=&e^{\pi i\alpha}e^{N/2}\frac{2\pi i}{e^{2\pi i\alpha}e^N-\id},\\
\Gamma(\id+\frac{N}{2\pi i})
\Gamma(\id-\frac{N}{2\pi i})
&=&e^{N/2}\frac{N}{e^N-\id}.
\end{eqnarray*}
Now calculate for $a\in H_\lambda,b\in H_{\oooo{\lambda}}$
with $\lambda=e^{-2\pi i \alpha}$ and $0<\alpha<1$
\begin{eqnarray*}
L^{nor}(Ga,e^{-N/2}Gb)
&=& L^{nor}(\Gamma(\alpha\id-\frac{N}{2\pi i})a,
e^{-N/2}\Gamma((1-\alpha)\id-\frac{N}{2\pi i})b)\\
&=& L^{nor}(a,e^{-N/2}\Gamma(\alpha\id+\frac{N}{2\pi i})
\Gamma(\id-(\alpha\id+\frac{N}{2\pi i}))b)\\
&=& L^{nor}(a,e^{\pi i\alpha}\frac{2\pi i}{M-\id}b)\\
&=& e^{\pi i \alpha}\cdot 2\pi i\cdot (-1)S(a,b).
\end{eqnarray*}
And calculate for $a,b\in H_1$
\begin{eqnarray*}
L^{nor}(Ga,e^{-N/2}Gb)
&=& L^{nor}(\Gamma(\id-\frac{N}{2\pi i})a,
e^{-N/2}\Gamma(\id-\frac{N}{2\pi i})b)\\
&=& L^{nor}(a,e^{-N/2}\Gamma(\id+\frac{N}{2\pi i})
\Gamma(\id-\frac{N}{2\pi i})b)\\
&=& L^{nor}(a,\frac{N}{M-\id}b)=S(a,b).
\end{eqnarray*}

(b) Definition \ref{t3.3} (c)(iii)\&(d) says here
\begin{eqnarray*}
S(F^pH_\lambda,F^{m+\theta(\lambda)+1-p}H_{\oooo{\lambda}})=0
\end{eqnarray*}
(recall $\theta(\lambda)$ in the notation \ref{t3.1}).
Part (b) follows from this, from part (a)
and from lemma \ref{t4.1} (a).

\medskip
(c) The spectral number $\alpha$ and the number
$\beta\in(0,1]$ with $e^{-2\pi i\beta}=\lambda$
satisfy $\alpha=m-p-1+\beta$. 
The positivity condition in definition \ref{t3.3} (c)(iv)
($\beta$) 
says for $a\in (I^{pq}_0)_\lambda-\{0\}$ as well as
for $a\in G((I^{pq}_0)_\lambda)-\{0\}$
\begin{eqnarray*}
0&<&i^{2p-m-\theta(\lambda)-l}\cdot S(a,N^l\oooo a).
\end{eqnarray*}
In the case $\lambda\neq 1$ this is because of \eqref{4.7}
\begin{eqnarray*}
&&i^{2p-m-l}\cdot \frac{-1}{2\pi i}\cdot e^{-\pi i\beta}\cdot
L^{nor}(Ga,e^{-N/2}G N^l\oooo{a})\\
&=& \frac{\Gamma(\beta)\Gamma(1-\beta)}{2\pi}
e^{\pi i(p-\frac{1}{2}(m+l-1))}e^{-\pi i(\alpha-m+p+1)}
L^{nor}(a,N^l\oooo{a})\\
&=& \frac{\Gamma(\beta)\Gamma(1-\beta)}{2\pi}
e^{\frac{1}{2}\pi i(m-l-1-2\alpha)}L^{nor}(a,N^l\oooo{a}).
\end{eqnarray*}
In the case $\lambda=1$ it is because of \eqref{4.8}
\begin{eqnarray*}
&&i^{2p-m-1-l}\cdot
L^{nor}(Ga,e^{-N/2}G N^l\oooo{a})\\
&=& e^{\pi i(p-\frac{1}{2}(m+1+l))}
L^{nor}(a,N^l\oooo{a})\\
&=& e^{\frac{1}{2}\pi i(m-l-1-2\alpha)}L^{nor}(a,N^l\oooo{a}).
\qquad\Box
\end{eqnarray*}

\bigskip
The Hodge filtration $F^\bullet$ is self-isotropic
with respect to $S$ by definition \ref{t3.3} (c)(iii)\&(d).
The Hodge filtration $G(F^\bullet)$ is self-isotropic 
with respect to the pairings $L^{sym}_{\kappa}$ in
definition \ref{t4.1} (a) by theorem \ref{t4.3} (b).

In section \ref{s6} we will see that $G(F^\bullet)$
behaves well with respect to a Thom-Sebastiani formula
in the singularity case. Below in theorem \ref{t4.6}
a special case is formalized and gives a certain
square root of a Tate twist for Steenbrink PMHS.

If $N\neq 0$ then the automorphism $G$ in definition
\ref{t4.2} (a) does not respect $H_\R$.
Then Deligne's $I^{pq}(G(F^\bullet))$ 
for $G(F^\bullet)$ are not equal
to the images $G(I^{pq}(F^\bullet))$ under $G$ of
Deligne's $I^{pq}(F^\bullet)$ for $F^\bullet$.
In view of \eqref{3.9}, the images $G(I^{pq}(F^\bullet))$
satisfy isotropy conditions for the pairings 
$L^{sym}_\kappa$, probably contrary to the 
$I^{pq}(G(F^\bullet))$. Therefore we worked in theorem
\ref{t4.3} (c) with the $I^{pq}(F^\bullet)$
and the images $G(I^{pq}(F^\bullet))$.

The next theorem \ref{t4.4} gives the classification of
the irreducible pieces in the Seifert form pair
$(H_\R,L^{nor})$ of a Steenbrink PMHS.
The proof uses theorem \ref{t3.8} and theorem \ref{t4.3}.
Part (b) recovers Nemethi's result \cite{Ne95}
that the isomorphism class of the Seifert form pair
$(H_\R,L^{nor})$ together with the number $m$ on one side and 
the spectral pairs modulo $2\Z$ in the first entry,
i.e. the data $\Spp\modd 2\Z\times\{0\}$, 
on the other side, determine each other.

\begin{theorem}\label{t4.4}
Let $(H_\R,H_\C,F^\bullet,W_\bullet,M,S)$ be a Steenbrink PMHS
of weight $m$.
Recall that $2\alpha+l+1-m$ is the distance from a
spectral pair ladder to its partner. 
It is an integer if and only if
$\lambda\in\{\pm 1\}$.

(a) For each ordered pair of spectral pair ladders or a single
spectral pair ladder the first spectral pair
is called $(\alpha,m+l)$, and then $\lambda:=e^{-2\pi i \alpha}$.
The Seifert form pair $(H_\R,L^{nor})$ (from definition 
\ref{t4.2} (c)) decomposes as follows.

It contains for each ordered pair of spectral pair ladders with
$\lambda\in S^1-\{0\}$ a Seifert form pair of the type
\begin{eqnarray}\label{4.12}
\Seif((-1)^{m+1}\lambda,2,l+1,e^{\frac{1}{2}\pi i(2\alpha+l+1-m)}).
\end{eqnarray}
It contains for each
pair of spectral pair ladders with odd distance
an irreducible Seifert form pair of the type
\begin{eqnarray}\label{4.13}
\Seif((-1)^{m+1}\lambda,2,l+1),\\
\textup{and then }(-1)^{m+1}\lambda=(-1)^{l+1}.\nonumber
\end{eqnarray}
It contains for each pair of spectral pair ladders
with even distance and each single spectral pair ladder 
(then the distance $2\alpha+l+1-m$ is 0) 
two respectively one Seifert form pair(s) of the type 
\begin{eqnarray}\label{4.14}
\Seif((-1)^{m+1}\lambda,1,l+1,
(-1)^{\frac{1}{2}(2\alpha+l+1-m)}),\\
\textup{and then }(-1)^{m+1}\lambda=(-1)^{l}.\nonumber
\end{eqnarray}

(b) $\Spp\modd 2\Z\times\{0\}$ and the isomorphism class
of $(H_\R,L^{nor})$ together with $m$ determine one another.
\end{theorem}

{\bf Proof:}
(a) By remark \ref{t3.7} we can suppose as in theorem
\ref{t3.8} that the Steenbrink PMHS is split.
Then we can consider the isometric triples in theorem \ref{t3.8}
and the corresponding Seifert form pairs with $L^{nor}$.

Remark that the monodromy of $L^{nor}$ is 
$M^{nor}=(-1)^{m+1}M$, so the eigenvalues change
from $\lambda$ to $(-1)^{m+1}\lambda$. 

If $\lambda\neq\pm 1$, the isometric triple and the
corresponding Seifert form pair are both irreducible.
Then \eqref{4.11} and theorem \ref{t2.9} give \eqref{4.12}.

In the other cases $\lambda\in\{\pm 1\}$.
Then by lemma \ref{t2.7}, the isometric triple
in \eqref{3.21} is irreducible if and only if
\begin{eqnarray*}
(l+1)+[m+\theta(\lambda)]_2+1\equiv 1(2),
\quad\textup{ i.e. }l+m+\theta(\lambda)\equiv 1(2).
\end{eqnarray*}
But
\begin{eqnarray*}
2\alpha&\equiv& \theta(\lambda)+1\modd 2\quad\textup{and then}\\
2\alpha+l+1-m&\equiv& l+m+\theta(\lambda)\modd 2.
\end{eqnarray*}
So, in the case of an odd distance $2\alpha+l+1-m$, 
the isometric triple in \eqref{3.21} is irreducible. 
By the proof of theorem \ref{t2.9} then also the corresponding 
Seifert form pair is irreducible. This gives \eqref{4.13}.

In the case of an even distance $2\alpha+l+1-m$,
the isometric triple and the Seifert form pair are
both reducible. 
Each pair of spectral pair ladders and each single spectral
pair ladder give two respectively one Seifert form pair
$\Seif((-1)^{m+1}\lambda,1,l+1,\varepsilon)$.
Here $\varepsilon=(-1)^{\frac{1}{2}(2\alpha+l+1-m)}$
because of \eqref{4.11} and theorem \ref{t2.9}.
This shows \eqref{4.14}.

\medskip
(d) It is rather easy to see that $\Spp\modd 2\Z\times\{0\}$
is equivalent to the spectral pair ladders modulo
$2\Z\times\{0\}$. Part (a) shows that these are equivalent
to the union of the corresponding Seifert form pairs
in $(H_\R,L)$ together with the number $m$.
The number $m$ is used to fix the symmetry point
$(\frac{m-1}{2},m)$ of $\Spp$.
\hfill$\Box$

\begin{remark}\label{t4.5}
In the case $N=0$ a Steenbrink PMHS can also be called a 
(pure) Steenbrink PHS. Then $G(F^\bullet)=F^\bullet$
and $I^{pq}_0=I^{pq}=H^{pq}$ with $q=m+\theta(\lambda)-p$.

Then the $\nu$ in \eqref{4.5} is $-\id$ on $H_1$.
Define $L^{nor}$ and $M^{nor}$ as in definition \ref{t4.2}.
$M^{nor}$ has on $H_\lambda$ the eigenvalue $\lambda\cdot
(-1)^{m+1}$.

Then on $H_\lambda$ the hermitian form $L^{herm}_\kappa$
from lemma \ref{t4.1} (b)
for any (of the two) $\kappa$ with 
$\kappa^2=\lambda\cdot (-1)^{m+1}$
is up to a constant equal to the hermitian form
$i^{-m-\theta(\lambda)} S(.,\oooo{.})$.

The Hodge decomposition $\bigoplus_p (H^{pq})_\lambda$
is then orthogonal with respect to $L^{herm}_\kappa$.
The positivity condition in \eqref{4.11} can then be written as
\begin{eqnarray}\label{4.15}
L^{herm}_{\exp(-\pi i(\alpha-\frac{m-1}{2}))}(a,a)>0
\end{eqnarray}
for $a\in H^{p,m+\theta(\lambda)-p}_\lambda-\{0\}$
and the spectral number $\alpha$ with
$e^{-2\pi i\alpha}=\lambda$ and $\lfloor m-\alpha\rfloor =p$.
\end{remark}

The following theorem \ref{t4.6} constructs from
a Steenbrink PMHS of weight $m$ a Steenbrink PMHS
of weight $m+1$, with the same underlying normalized
Seifert form pair $(H_\R,L^{nor})$.
In the singularity case it corresponds to a suspension:
one goes from a singularity $f(x_0,...,x_m)$ to a
singularity $f(x_0,...,x_m)+x_{m+1}^2$, 
see remark \ref{t6.6} (iii).
It can be seen as a square root of a Tate twist.

\begin{theorem}\label{t4.6}
Let $(H_\R,H_\C,F^\bullet,W_\bullet,M,S)$ be a Steenbrink
PMHS of weight $m$. We construct a Steenbrink PMHS
$(\www H_\R,\www H_\C,\www F^\bullet, \www W_\bullet,
\www M,\www S)$ of weight $\www m=m+1$ as follows:
\begin{eqnarray*}
\www H_\R=H_\R,\quad \www H_\C=H_\C, \quad\www M=-M,
\quad \www M_s=-M_s,\quad \www N=N,\\
\www H_\lambda=H_{-\lambda},\quad \www W_\bullet\www H_{\neq 1}
=W^{(N,m+1)}_\bullet \www H_{\neq 1},\quad
\www W_\bullet \www H_1=W^{(N,m+2)}_\bullet \www H_1,
\end{eqnarray*}
\begin{eqnarray*}
\www\nu &:=&\left\{\begin{array}{ll}
\frac{1}{\www M-\id}&\textup{on }\www H_{\neq 1}
\quad(\textup{eigenvalues }\neq 1\textup{ w.r.t. }\www M)\\
\frac{-N}{\www M-\id}&\textup{on }\www H_1
\quad(\textup{eigenvalue }1\textup{ w.r.t. }\www M),
\end{array}\right.\\
\www S(a,b)&:=& -L^{nor}(a,\www\nu b),
\end{eqnarray*}
\begin{eqnarray*}
F^p\www H_\lambda &:=& (G^{(\alpha+\frac{1}{2})})^{-1}
G^{(\alpha)} F^{p+1}\www H_\lambda
\qquad\textup{if }-\lambda=e^{-2\pi i \alpha},
0<\alpha\leq\frac{1}{2},\\
F^p\www H_\lambda &:=& (G^{(\alpha-\frac{1}{2})})^{-1}
G^{(\alpha)} F^{p}\www H_\lambda
\qquad\textup{if }-\lambda=e^{-2\pi i \alpha},
\frac{1}{2}<\alpha\leq 1.
\end{eqnarray*}
Then $\www\Spp=\Spp+(\frac{1}{2},1)$. 
The two Steenbrink PMHS have the same underlying
Seifert form pair $(H_\R,L^{nor})$.
Carrying out the construction twice, leaves
$(H_\R,H_\C,M,S)$ invariant and gives
$\www{\www m}=m+2$, 
$\www{\www W}_\bullet=W_{\bullet+2}$, 
$\www{\www F^\bullet} =F^{\bullet+1}$,
so it is a Tate twist.
\end{theorem}

{\bf Proof:}
The proof uses theorem \ref{t4.3}. 
We leave the details to the reader. Compare also 
corollary \ref{t6.5} and remark \ref{t6.6} (iii).
\hfill$\Box$

\section{Fourier-Laplace transformation and pairings
on a bundle on $\C$ with regular singular connection
on $(\C,0)$}\label{s5}
\setcounter{equation}{0}

\noindent 
We will present an equivalence between three types of pairings
and additional data, a polarizing form plus a monodromy,
a Seifert form, and a pairing on a flat bundle
on $\C^*$. Then we will consider holomorphic sections
with moderate growth in the bundle and study a 
Fourier-Laplace transformation on them.

This will make the meaning of the automorphism $G$ in
definition \ref{t4.2} transparent. 
Theorem \ref{t5.2} will also fill the equivalence
with life, by nice formulas which connect the pairings.
Theorem \ref{t5.2} was stated in \cite{He03} as proposition 7.7,
but the proof was essentially omitted.

\begin{lemma}\label{t5.1}
The following three data are equivalent.

$(\alpha)$ $(H_\R,M,S,m)$. Here $H_\R$ is a finite dimensional
$\R$-vector space, $M$ is an automorphism on it with eigenvalues
in $S^1$, called {\rm monodromy}.
$m\in\Z$. And $S$ is an $M$-invariant nondegenerate bilinear form.
On $H_{\neq 1}$ it is $(-1)^m$-symmetric.
On $H_1$ it is $(-1)^{m+1}$-symmetric.

$(\beta)$. $(H_\R,L,m)$. Here $(H_\R,L)$ is a Seifert form 
pair such that the eigenvalues of $L$ are in $S^1$,
and $m\in\Z$.

$(\gamma)$ $(H^{bun}_\R\to\C^*,\nnn,P,m)$.
Here $H^{bun}_\R\to\C^*$ is a bundle of $\R$-vector spaces
on $\C^*$ with flat connection $\nnn$, whose monodromy
has eigenvalues in $S^1$. Its
complexification is a holomorphic flat bundle and is
denoted $H^{bun}_\C\to\C^*$.
Again $m\in\Z$. And $P$ is a flat and nondegenerate and
$(-1)^{m+1}$-symmetric pairing
\begin{eqnarray}\label{5.1}
P:H^{bun}_{\R,z}\times H^{bun}_{\R,-z}\to i^{m+1}\cdot\R
\quad\textup{ for }z\in \C^*.
\end{eqnarray}

\medskip
From $(\alpha)$ to $(\beta)$: $L:=L^{nor}$ in \eqref{4.6},
using \eqref{4.5}.

From $(\beta)$ to $(\gamma)$: Define a flat bundle
$H^{bun}_\R\to\C^*$ with monodromy $(-1)^{m+1}M$.
Then $L:H^{bun}_{\R,z}\times H^{bun}_{\R,z}\to \R$
is defined on each fiber and is flat.
Define $P$ by 
\begin{eqnarray}\label{5.2}
P(a,b):= \frac{1}{(2\pi i)^{m+1}}\cdot L(a,\gamma_{-\pi}b).
\end{eqnarray}
Here $\gamma_{-\pi}:H^{bun}_{\R,z}\to H^{bun}_{\R,-z}$
is the isomorphism by flat shift in mathematically
negative direction.

One goes from $(\gamma)$ to $(\beta)$ and from $(\beta)$
to $(\alpha)$ by inverting these constructions.
\end{lemma}

{\bf Proof:}
Lemma \ref{t2.3} and lemma \ref{t2.4} show the equivalence
of $(\alpha)$ and $(\beta)$. 

From $(\beta)$ to $(\gamma)$: $P$ is well defined and
nondegenerate and flat because $L$ has these properties.
It is $(-1)^{m+1}$-symmetric because of
$L(b,a)=L(Ma,b)$: For $a\in H^{bun}_{\R,z}$, 
$b\in H^{bun}_{\R,-z}$
\begin{eqnarray}
&& (2\pi i)^{m+1}\cdot P(b,a)
=L(b,\gamma_{-\pi}a)=L(M\gamma_{-\pi}a,b) \nonumber\\
&=& (-1)^{m+1}L((-1)^{m+1}M\gamma_{-\pi}a,b)
=(-1)^{m+1}L(\gamma_\pi a,b)  \nonumber\\
&=& (-1)^{m+1}L(a,\gamma_{-\pi}b)
=(2\pi i)^{m+1}\cdot (-1)^{m+1}\cdot P(a,b).\label{5.3}
\end{eqnarray}
Here $\gamma_{\pi}:H^{bun}_{\R,z}\to H^{bun}_{\R,-z}$
is the isomorphism by flat shift in mathematically
positive direction.

From $(\gamma)$ to $(\beta)$: Define $L$ on any fiber of
$H^{bun}$ by $L(a,b):=(2\pi i)^{m+1}\cdot P(a,\gamma_\pi b)$.
The $(-1)^{m+1}$-symmetry of $P$ gives
$L(Ma,b)=L(b,a)$, by inverting the calculation \eqref{5.3}.
Take $H_\R:=H^{bun}_{\R,z}$ for an arbitrary $z\in\C^*$.
\hfill$\Box$ 

\bigskip
Now we consider all the data in lemma \ref{t5.1}.
Before we come to theorem \ref{t5.2}, we have to
describe the {\it elementary sections} and the
{\it Kashiwara-Malgrange $V$-filtration}.
Of course, this is standard and can be found at many places,
e.g. \cite{SS85}\cite{AGV88}\cite{He02}.

The space of flat multivalued global sections in
$H^{bun}_\C\to\C^*$ is denoted $H^\infty_\C$.
It can be identified in a non-unique way with $H_\C$.
It comes with a monodromy which is then identified with
$(-1)^{m+1}M$. Now $H^\infty_\lambda$ means the
generalized eigenspace with respect to this monodromy.
$H^\infty$ also comes with a real subspace $H^\infty_\R$.

Any global flat multivalued section $A\in H^\infty_\lambda$
and any choice of $\alpha\in\R$ with $e^{-2\pi i \alpha}=\lambda$
leads to a holomorphic univalued section with specific 
growth condition at $0\in\Delta$, the 
{\it elementary section} $es(A,\alpha)$ with
\begin{eqnarray*}
es(A,\alpha)(\tau):=
e^{\log\tau(\alpha-\frac{N}{2\pi i})}\cdot A(\log \tau).
\end{eqnarray*}
Denote by $C^\alpha$ the $\C$-vector space of all
elementary sections with fixed $\alpha$ and $\lambda$.
The map $es(.,\alpha):H^\infty_\lambda\to C^\alpha$
is an isomorphism. The space 
$V^{mod}:=\bigoplus_{\alpha\in(-1,0]}\C\{\tau\}[\tau^{-1}]
\cdot C^\alpha$ is the space of all germs at 0 of the sheaf
of holomorphic sections on the flat cohomology bundle with 
moderate growth at 0. The Kashiwara-Malgrange $V$-filtration
is given by the subspaces
\begin{eqnarray*}
V^\alpha:=\bigoplus_{\beta\in[\alpha,\alpha+1)}\C\{\tau\}
\cdot C^\beta,\quad
V^{>\alpha}:=\bigoplus_{\beta\in(\alpha,\alpha+1]}\C\{\tau\}
\cdot C^\beta.
\end{eqnarray*}
It is a decreasing filtration by free $\C\{\tau\}$-modules
of rank $\mu$ with 
$\Gr_V^\alpha=V^\alpha/V^{>\alpha}\cong C^\alpha$. And
\begin{eqnarray*}
\tau:C^\alpha\to C^{\alpha+1}\textup{ bijective},&&
\tau\cdot es(A,\alpha)=es(A,\alpha+1),\\
\paa_{\tau}:C^\alpha\to C^{\alpha-1}\textup{ bijective}
&&\textup{if }\alpha\neq 0,\\
\tau\paa_\tau-\alpha:C^\alpha\to C^\alpha\textup{ nilpotent},&&
(\tau\paa_\tau-\alpha)es(A,\alpha)=es(\frac{-N}{2\pi i}A,
\alpha).
\end{eqnarray*}

\begin{theorem}\label{t5.2}
\cite[Proposition 7.7]{He03}

(a) Let $\tau$ and $z$ both be coordinates on $\C$. For 
$\alpha>0$ and $A\in H^\infty_{e^{-2\pi i\alpha}}$,
the Fourier-Laplace transformation $FL$ with
\begin{eqnarray}\label{5.4}
FL(es(A,\alpha-1)(\tau))(z):=
\int_0^{\infty\cdot z}e^{-\tau/z}\cdot es(A,\alpha-1)(\tau)d\tau
\end{eqnarray}
is well defined and maps the elementary section 
$es(A,\alpha-1)(\tau)$ in $\tau$ to the elementary section
\begin{eqnarray}\label{5.5}
FL(es(A,\alpha-1)(\tau))(z)=es(G^{(\alpha)}A,\alpha)(z)
\end{eqnarray}
in $z$.

(b) For $0<\alpha<1$ and $A\in H^\infty_{e^{-2\pi i\alpha}}$,
$B\in H^\infty_{e^{2\pi i\alpha}}$,
\begin{eqnarray}\label{5.6}
&& P(es(G^{(\alpha)}A,\alpha)(z),es(G^{(1-\alpha)}B,1-\alpha)(-z)\\
&=& \frac{z}{(2\pi i)^{m+1}}\cdot e^{\pi i(1-\alpha)}
\cdot L^{nor}(G^{(\alpha)}A(\log z),e^{-N/2}G^{(1-\alpha)}
B(\log z))\nonumber \\
&=& \frac{z}{(2\pi i)^{m}}\cdot S(A,b).\nonumber 
\end{eqnarray}
For $A,B\in H^\infty_1$, 
\begin{eqnarray}\label{5.7}
&& P(es(G^{(1)}A,1)(z),es(G^{(1)}B,1)(-z)\\
&=& \frac{-z^2}{(2\pi i)^{m+1}}
\cdot L^{nor}(G^{(1)}A(\log z),e^{-N/2}G^{(1)}
B(\log z))\nonumber \\
&=& \frac{-z^2}{(2\pi i)^{m+1}}\cdot S(A,b).\nonumber 
\end{eqnarray}
\end{theorem}

{\bf Proof:}
As the proof was not carried out in \cite{He03},
we give it here.

(a) The Gamma function satisfies for $\alpha>0$ the identity
\begin{eqnarray*}
\left(\frac{d}{d\alpha}\right)^k
\left(\Gamma(\alpha)z^\alpha\right)
=\int_0^{\infty\cdot z}e^{-\tau/z}\cdot 
\tau^{\alpha-1}(\log\tau)^k d\tau.
\end{eqnarray*}
In the next calculation $j+l=k$,
\begin{eqnarray*}
&& FL(es(A,\alpha-1)(\tau))(z)\\
&=& \int_0^{\infty\cdot z}e^{-\tau/z}\cdot \tau^{\alpha-1}
\cdot\sum_{k\geq 0}\frac{1}{k!}(\log\tau)^k
\left(\frac{-N}{2\pi i}\right)^k A(\log\tau)d\tau\\
&=& \sum_{k\geq 0}\frac{1}{k!}\left(\frac{-N}{2\pi i}\right)^k
A(\log z)\cdot \left(\frac{d}{d\alpha}\right)^k
\left(\Gamma(\alpha)z^\alpha\right)\\
&=& \sum_{k\geq 0}\frac{1}{k!}\left(\frac{-N}{2\pi i}\right)^k
A(\log z)\cdot 
\sum_{l=0}^k \begin{pmatrix}k\\l\end{pmatrix}
\Gamma^{(l)}(\alpha)\cdot (\log z)^{k-l}\cdot z^\alpha\\
&=& \sum_{j,l\geq 0}z^\alpha\frac{1}{j!}
\left(\log z\cdot\frac{-N}{2\pi i}\right)^j\cdot 
\frac{1}{l!}\Gamma^{(l)}(\alpha)\left(\frac{-N}{2\pi i}\right)^l
A(\log z)\\
&=& e^{\log z(\alpha-\frac{N}{2\pi i})}
\Gamma(\alpha-\frac{N}{2\pi i}) A(\log z)
= es(G^{(\alpha)}A,\alpha)(z).
\end{eqnarray*}

(b) 
For $0<\alpha<1$ and $A\in H^\infty_{e^{-2\pi i\alpha}}$,
$B\in H^\infty_{e^{2\pi i\alpha}}$
\begin{eqnarray*}
&&(2\pi i)^{m+1}\cdot P(es(A,\alpha)(z),es(B,1-\alpha)(-z))\\
&=&
L^{nor}(es(A,\alpha)(z),\gamma_{-\pi}es(B,1-\alpha)
(e^{\pi i}z))\\
&=& 
L^{nor}(e^{\log z(\alpha-\frac{N}{2\pi i})}A(\log z),
\gamma_{-\pi}e^{(\pi i+\log z)(1-\alpha-\frac{N}{2\pi i})}
B(\pi i+\log z))\\
&=&
L^{nor}(A(\log z),e^{\log z(\alpha+\frac{N}{2\pi i})}
e^{(\pi i+\log z)(1-\alpha-\frac{N}{2\pi i})}B(\log z))\\
&=& 
z\cdot e^{\pi i(1-\alpha)}\cdot
L^{nor}(A(\log z),e^{-N/2}B(\log z)).
\end{eqnarray*}
For $A,B\in H^\infty_1$,
\begin{eqnarray*}
&&(2\pi i)^{m+1}\cdot P(es(A,1)(z),es(B,1)(-z))\\
&=& 
L^{nor}(es(A,1)(z),\gamma_{-\pi}es(B,1)(e^{\pi i}z))\\
&=& 
L^{nor}(e^{\log z(1-\frac{N}{2\pi i})}A(\log z),
\gamma_{-\pi}e^{(\pi i+\log z)(1-\frac{N}{2\pi i})}
B(\pi i+\log z))\\
&=& 
L^{nor}(A(\log z),e^{\log z(1+\frac{N}{2\pi i})}
e^{(\pi i+\log z)(1-\frac{N}{2\pi i})}B(\log z))\\
&=& 
(-z^2)\cdot
L^{nor}(A(\log z),e^{-N/2}B(\log z)).
\end{eqnarray*}
The equalities involving $S$ follow now with theorem 
\ref{t4.3} (a). \hfill$\Box$

\section{Isolated hypersurface singularities}\label{s6}
\setcounter{equation}{0}

\noindent
Our main motivation for this paper is the study of 
isolated hypersurface singularities. Each comes 
with its Milnor lattice, a $\Z$-lattice with an integer
valued Seifert form.

A singularity comes also with a {\it signed} Steenbrink PMHS
and thus with a polarizing form and spectral pairs. 
This section defines and
names all these data and states results.
Most of it is well known. The main new point is the
correction in corollary \ref{t6.5} of a Thom-Sebastiani formula
for the signed Steenbrink PMHS in \cite[ch. 8]{SS85}.
But this correction requires the material
in the sections \ref{s4} and \ref{s5}.

An {\it isolated hypersurface singularity}
(short: {\it singularity})
is a holomorphic function germ $f:(\C^{m+1},0)\to (\C,0)$ 
with an isolated singularity at $0$. Its {\it Milnor number} 
$$\mu:=\dim\OO_{\C^{m+1},0}/(\frac{\paa f}{\paa x_0},...,
\frac{\paa f}{\paa x_m})$$ 
is finite. 
For the following notions and facts compare \cite{AGV88} and 
\cite{Eb07}.
A {\it good representative} of $f$ has to be defined with some 
care \cite{Mi68}\cite{AGV88}\cite{Eb07}. It is $f:X\to \Delta$
with $X=\{x\in\C^{m+1}\, |\, |x|<\varepsilon\}\cap f^{-1}(\Delta)$ 
for a sufficiently small $\varepsilon>0$ and
$\Delta = \{\tau\in\C\, |\, |\tau|<\delta\}$ 
a small disk around 0 (first choose $\varepsilon$, then
$\delta$). 
Then $f:X'\to \Delta'$ with $X'=X-f^{-1}(0)$ and 
$\Delta'=\Delta-\{0\}$ is a locally trivial $C^\infty$-fibration,
the  {\it Milnor fibration}. Each fiber has the
homotopy type of a bouquet of $\mu$ $n$-spheres \cite{Mi68}.

Therefore the (for $m=0$ reduced) 
middle homology groups are {}\\{}
$H_m^{(red)}(f^{-1}(\tau),\Z) \cong \Z^\mu$ for $\tau\in \Delta'$.
Each comes equipped with an intersection form $I$, 
which is a datum of one fiber,
a monodromy $M$ and a Seifert form $L$, which come from the 
Milnor fibration,
see \cite[I.2.3]{AGV88} for their definitions.
$M$ is a quasiunipotent automorphism, $I$ and $L$ are 
bilinear forms with values in $\Z$,
$I$ is $(-1)^m$-symmetric, and $L$ is unimodular. $
L$ determines $M$ and $I$ because of the formulas
\cite[I.2.3]{AGV88}
\begin{eqnarray}\label{6.1}
L(Ma,b)&=&(-1)^{m+1}L(b,a),\\ \label{6.2}
I(a,b)&=&-L(a,b)+(-1)^{m+1}L(b,a)=L((M-\id)a,b).
\end{eqnarray}
The Milnor lattices $H_m(f^{-1}(\tau),\Z)$ for all 
Milnor fibrations $f:X'\to \Delta'$ and then all 
$\tau\in\R_{>0}\cap T'$ are canonically isomorphic,
and the isomorphisms respect $M$, $I$ and $L$. 
This follows from Lemma 2.2 in \cite{LR73}. 
These lattices are identified and called 
{\it Milnor lattice} $Ml(f)$.

A result of Thom and Sebastiani
compares the Milnor lattices and monodromies of 
the singularities $f=f(x_0,...,x_m),g=g(y_0,...,y_n)$ and
$f+g=f(x_0,...,x_m)+g(x_{m+1},...,x_{m+n+1})$.
There is an extension by Deligne for the Seifert form
\cite[I.2.7]{AGV88}. It is restated here.
There is a canonical isomorphism
\begin{eqnarray}\label{6.3}
\Phi:Ml(f+g)&\stackrel{\cong}{\longrightarrow} 
&Ml(f)\otimes Ml(g),\\
\textup{with } M(f+g)&\cong & M(f)\otimes M(g) 
\label{6.4}\\
\textup{and } 
L(f+g)&\cong& (-1)^{(m+1)(n+1)}\cdot L(f)\otimes L(g).
\label{6.5}
\end{eqnarray}
This motivates the definition of the {\it normalized
Seifert form}  and the {\it normalized monodromy} 
on the Milnor lattice $Ml(f)$
\begin{eqnarray}\label{6.6}
L^{hnor}(f)&:=&(-1)^{(m+1)(m+2)/2}\cdot L(f),\\
M^{hnor}(f)&:=& (-1)^{m+1}\cdot M(f)\label{6.7}
\end{eqnarray}
because then 
\begin{eqnarray}\label{6.8}
L^{hnor}(f+g)&\cong& L^{hnor}(f)\otimes L^{hnor}(g),\\
M^{hnor}(f+g)&\cong& M^{hnor}(f)\otimes M^{hnor}(g)
\label{6.9}
\end{eqnarray}
and $M^{hnor}$ is the monodromy of $L^{hnor}$ in the sense of 
lemma \ref{t2.3} (a). 

In the special case $g=x_{m+1}^2$,
the function germ 
$f+g=f(x_0,...,x_m)+x_{m+1}^2\in \OO_{\C^{m+2},0}$
is called {\it stabilization} or {\it suspension} of $f$. 
As there are only two isomorphisms $Ml(x_{m+1}^2)\to\Z$, 
and they differ by a sign, there are two equally canonical
isomorphisms $Ml(f)\to Ml(f+x_{m+1}^2)$,
and they differ just by a sign. 
Therefore automorphisms and bilinear forms on $Ml(f)$ 
can be identified with automorphisms and bilinear forms on 
$Ml(f+x_{m+1}^2)$. In this sense \cite[I.2.7]{AGV88}
\begin{eqnarray}\label{6.10}
L(f+x_{m+1}^2) &=& (-1)^m\cdot L(f),\\ 
M(f+x_{m+1}^2)&=& - M(f),\label{6.11}\\
L^{hnor}(f+x_{m+1}^2) &=&  L^{hnor}(f),\label{6.12}\\
M^{hnor}(f+x_{m+1}^2)&=& M^{hnor}(f).\label{6.13}
\end{eqnarray}

Denote by $H^\infty_\C$ the $\mu$-dimensional vector space
of global flat multivalued sections in the flat cohomology
bundle $\bigcup_{\tau\in\Delta'}H^m(f^{-1}(\tau),\C)$
(reduced cohomology for $m=0$). 
It comes equipped with a $\Z$-lattice $H^\infty_\Z$
and a real subspace $H^\infty_\R$ and a monodromy which
is also denoted by $M$. 

There is a natural {\it signed} Steenbrink PMHS
({\it signed}: definition \ref{t6.1} below)
$(H^\infty_\C,H^\infty_\R,F^\bullet_{St},W_\bullet,M,S)$
on $H^\infty_\C$. The weight filtration is
$W^{(N,m)}_\bullet H^\infty_{\neq 1}$ on $H^\infty_{\neq 1}$
and $W^{(N,m+1)}_\bullet H^\infty_1$ on $H^\infty_1$
(see lemma \ref{t3.2} (a) for $W_\bullet^{(N,m)}$).
The Hodge filtration was defined first by Steenbrink using
resolution of singularities \cite{St77}.
Then Varchenko \cite{Va80} constructed a closely related
Hodge filtration $F^\bullet_{Va}$ from the Brieskorn lattice
$H_0''(f)$ (definition below). Scherk and Steenbrink
\cite{SS85} and M. Saito \cite{SaM82} modified this 
construction to recover $F^\bullet_{St}$.
Below we explain the Brieskorn lattice and this modified
construction.

But first we give the polarizing form $S$. 
The lattice $H^\infty_\Z$
can be identified with the dual
$Ml(f)^\vee=\Hom(Ml(f),\Z)$ of the Milnor lattice $Ml(f)$,
and thus it comes equipped with the dual Seifert form
$L^\vee$ (using the notations \ref{t2.11})
of the Seifert form $L$ on $Ml(f)$. Define the normalized
Seifert form $L^{nor}$ on $H^\infty_\Z$ by
\begin{eqnarray}\label{6.14}
L^{nor}:=(-1)^{(m+1)(m+2)/2}L^\vee=(L^{hnor})^\vee,
\end{eqnarray}
an $M$-invariant automorphism $\nu:H^\infty_\Q\to H^\infty_\Q$
\begin{eqnarray}\label{6.15}
\nu :=\left\{\begin{array}{ll}
\frac{1}{M-\id}&\textup{on }H^\infty_{\neq 1},\\
\frac{-N}{M-\id}&\textup{on }H^\infty_1,\end{array}\right.
\end{eqnarray}
and the $M$-invariant polarizing form 
$S:H^\infty_\Q\times H^\infty_\Q\to\Q$ by 
\begin{eqnarray}\label{6.16}
S(a,b):=-L^{nor}(a,\nu b).
\end{eqnarray}
$L^{nor}$ and $S$ are related by the equivalence
in lemma \ref{t5.1}.
Therefore $S$ is $(-1)^m$-symmetric on 
$H^\infty_{\neq 1}$ and $(-1)^{m+1}$-symmetric on $H^\infty_1$.
The restriction to $H^\infty_{\neq 1}$ is
$(-1)^{m(m+1)/2}\cdot I^\vee$, where $I^\vee$ on 
$H^\infty_{\neq 1}$ is dual to $I$
(which is non-degenerate on $Ml(f)_{\neq 1}$).
This follows from \eqref{2.31} in corollary \ref{t2.12}.

Steenbrink had this restriction to $H^\infty_{\neq 1}$ of $S$,
but not the part on $H^\infty_1$. 
That part was defined with a sign mistake in \cite{He99}
and correctly in \cite{He02}. 
The same sign mistake led to the claim in \cite{He99}
that $(H^\infty_\R,H^\infty_\C,F^\bullet_{St},W_{\bullet},M,S)$
is a Steenbrink PMHS of weight $m$. 
But it is (as stated correctly in \cite{He02})
a {\it signed} Steenbrink PMHS of weight $m$.

\begin{definition}\label{t6.1}
A tuple $(H_\R,H_\C,F^\bullet,W_{\bullet},M,S)$
is a {\it signed Steenbrink PMHS} of weight $m$ if 
$(H_\R,H_\C,F^\bullet,W_{\bullet},M_s\cdot e^{-N},S)$
is a Steenbrink PMHS of weight $m$.
\end{definition}

\begin{remarks}\label{t6.2}
(i) The only difference between a Steenbrink PMHS and a signed
Steenbrink PMHS is that the positivity condition in 
definition \ref{t3.3} (c)(iv)($\beta$) 
(see the notations \ref{t3.1} for $\theta(\lambda)$)
\begin{eqnarray*}
i^{2p-m-\theta(\lambda)-l}\cdot S(a,N^l\oooo{a})>0
\end{eqnarray*}
for $a\in (F^pP_{m+l}\cap
\oooo{F^{m+\theta(\lambda)+l-p}P_{m+l}})_\lambda-\{0\}$
has to be replaced by the positivity condition
\begin{eqnarray*}
i^{2p-m-\theta(\lambda)-l}\cdot S(a,(-N)^l\oooo{a})>0.
\end{eqnarray*}
This changes the sign in the case of a Jordan block of 
even size, i.e. in the case of a pair of spectral
pair ladders (or a single one) of even length $l+1$.

\medskip
(ii) This leads to obvious variants of the theorems
\ref{t3.8}, \ref{t4.3} and \ref{t4.4} for signed
Steenbrink PMHS: 
In \eqref{3.21} and \eqref{3.22} the last entry
$\varepsilon\in\{\pm 1\}$ in the isometric triples
has to be replaced by $-\varepsilon$ if $l+1$ is even.
The factor in \eqref{4.11} and the last entry in the
Seifert form pairs in \eqref{4.12} and \eqref{4.14}
have to be multiplied by $-1$ if $l+1$ is even.

\medskip
(iii) We did not work from the beginning only with
signed Steenbrink PMHS because also Steenbrink PMHS
naturally appear.
An M-tame function on an affine manifold of dimension $m+1$ 
leads by work of Sabbah \cite{Sa98}\cite{NS99} to a 
Steenbrink MHS $(H_\R,H_\C,F^\bullet_{Sa},W_\bullet,M)$ 
of weight $m$. Then for $S$ defined as above, the tuple
$(H_\R,H_\C,G^{-1}F^\bullet_{Sa},W_\bullet,M,S)$
is a Steenbrink PMHS of weight $m$ \cite[theorem 7.3]{HS07}.

\medskip
(iv) A Steenbrink MHS of weight $m$ with polarizing form $S$
is a Steenbrink PMHS respectively a signed Steenbrink PMHS
if and only if $e^{zN}F^\bullet$ is for $\Imm z\gg 0$
respectively for $\Imm z\ll 0$ the Hodge filtration
of a sum of pure polarized Hodge structures 
(of weight $m$ on $H_{\neq 1}$ and of weight
$m+1$ on $H_1$) \cite{CKS86} (see also \cite{He03}\cite{HS07}).
This was lifted  in \cite{HS07} to Sabbah orbits respectively
nilpotent orbits of TERP-structures.
\end{remarks}

\begin{examples}\label{t6.3}
(i) The Laurent polynomial $x_0+\frac{1}{x_0}$ (so $m=0$) 
is an M-tame function on $\C^*$.
It is the mirror partner of the quantum cohomology of $\P^1$.
It has two $A_1$-singularities, so the global Milnor number
is $\mu=2$. Here the Milnor lattice has to be replaced by a
$\Z$-lattice of rank 2 of Lefschetz thimbles. This comes
equipped with a Seifert form.
For a suitable basis the matrix of the Seifert form is
$-\begin{pmatrix}1&0\\2&1\end{pmatrix}=:-S^t$.
Here the normalized Seifert form is $L^{hnor}=-L$, 
and the normalized monodromy is $M^{hnor}=-M$.
It has the matrix 
\begin{eqnarray*}
S^{-1}S^t
=\begin{pmatrix}-3&-2\\2&1\end{pmatrix}
=-\exp\begin{pmatrix}2&2\\-2&-2\end{pmatrix}
\end{eqnarray*}
with one $2\times 2$ Jordan block and eigenvalue $-1$
and matrix $\begin{pmatrix}2&2\\-2&-2\end{pmatrix}$
of its nilpotent part $N$. For the vector $a$ represented
by $\begin{pmatrix}1\\0\end{pmatrix}$ one finds
\begin{eqnarray}\label{6.17}
L^{hnor}(a,Na)=\begin{pmatrix}1&0\end{pmatrix}
\begin{pmatrix}1&0\\2&1\end{pmatrix}
\begin{pmatrix}2&2\\-2&-2\end{pmatrix}
\begin{pmatrix}1\\0\end{pmatrix}
=2>0.
\end{eqnarray}
Thus by theorem \ref{t2.9} the Seifert form pair
$(H_\R,L^{hnor})$ is of type $\Seif(-1,1,2,1)$.
This is in accordance with the fact that here we have a
Steenbrink PMHS of weight one and with \eqref{4.14},
which predict this type $\Seif(-1,1,2,1)$.

\medskip
(ii) Each hyperbolic surface (so $m=2$) 
singularity $T_{pqr}$ (with 
$\frac{1}{p}+\frac{1}{q}+\frac{1}{r}<1$)
has a rank 2 $\Z$-sublattice $Ml(f)_1\cap Ml(f)$ of the
Milnor lattice $Ml(f)$. For a suitable basis
$\uuuu{a}=(a_1,a_2)$ of this sublattice, the matrix of the
Seifert form $L$ is by \cite[(29)]{GH17}
\begin{eqnarray*}
L(\uuuu{a}^t,\uuuu{a})&=&
\begin{pmatrix}0&-\chi\\ \chi&\frac{\chi^2}{2}(\kappa-1)
\end{pmatrix} =:S^t\\
\textup{where }\kappa&:=&\frac{1}{p}+\frac{1}{q}+\frac{1}{r}<1,
\quad \chi:=\lcm(p,q,r).
\end{eqnarray*}
Here the normalized Seifert form is 
$L^{hnor}=L$ and is given by the matrix $S^t$. Its monodromy
is $M^{hnor}=-M$ and  has on $Ml(f)_1$ the matrix
\begin{eqnarray*}
S^{-1}S^t&=&
\frac{1}{\chi^2}
\begin{pmatrix}\frac{\chi^2}{2}(\kappa-1)&-\chi\\ \chi&0
\end{pmatrix}
\begin{pmatrix}0&-\chi\\ \chi&\frac{\chi^2}{2}(\kappa-1)
\end{pmatrix}\\
&=&-\exp \begin{pmatrix}0&\chi(\kappa-1)\\0&0\end{pmatrix}.
\end{eqnarray*}
Its nilpotent part $N$ on $Ml(f)_1$ has the matrix 
$\begin{pmatrix}0&\chi(\kappa-1)\\ 0&0\end{pmatrix}$.
For the vector $a$ represented
by $\begin{pmatrix}0\\1\end{pmatrix}$ one finds
\begin{eqnarray}
L^{hnor}(a,Na)&=&\begin{pmatrix}0&1\end{pmatrix}
\begin{pmatrix}0&-\chi\\ \chi&\frac{\chi^2}{2}(\kappa-1)
\end{pmatrix}
\begin{pmatrix}0&\chi(\kappa-1)\\0&0\end{pmatrix}
\begin{pmatrix}0\\1\end{pmatrix}\nonumber\\
&=&\chi^2(\kappa-1)<0.\label{6.18}
\end{eqnarray}
Thus by theorem \ref{t2.9} the Seifert form pair
$(Ml(f)_1\cap Ml(f)_\R,L^{hnor})$ is of type $\Seif(-1,1,2,-1)$.
This is in accordance with the fact that here we have a signed
Steenbrink PMHS of weight one and with the variant of \eqref{4.14},
which predict this type $\Seif(-1,1,2,-1)$.
\end{examples}

\bigskip
Now we apply the notations from section \ref{s5}
to the cohomology bundle
$\bigcup_{\tau\in\Delta'}H^m(f^{-1}(\tau),\C)$, 
i.e. the $V$-filtration
and the spaces $C^\alpha$ and the isomorphisms
$es(.,\alpha):H^\infty_{e^{-2\pi i\alpha}}\to C^\alpha$.

The {\it Brieskorn lattice} is a free $\C\{\tau\}$-module
$H_0''(f)\subset V^{>-1}$ which had first been defined
by Brieskorn \cite{Br70}. The name {\it Brieskorn lattice}
is due to \cite{SaM89}, the notation $H_0''$ is from \cite{Br70}.
The Brieskorn lattice is generated by germs of sections
$s[\omega]$ from holomorphic $(m+1)$-forms 
$\omega\in\Omega^{m+1}_X$: Integrating the 
Gelfand-Leray form $\frac{\omega}{df}|_{f^{-1}(\tau)}$
over cycles in $H_m(f^{-1}(\tau),\C)$ gives a holomorphic
section $s[\omega]$ in the cohomology bundle, whose germ
$s[\omega]_0$ at 0 is in fact in $V^{>-1}$
(this was proved first by Malgrange).

Steenbrink's Hodge filtration $F^\bullet_{St}H^\infty$
can be recovered as follows \cite{SS85}\cite{SaM82}.
Consider for $\lambda\in S^1$ the unique
$\alpha\in (0,1]$ with $\lambda=e^{-2\pi i\alpha}$. Then
\begin{eqnarray}\label{6.19}
F^p_{St}H^\infty_\lambda = es(.,\alpha-1)(\tau)^{-1}
\Bigl(\paa_\tau^{m-p}\Gr_V^{m-p+\alpha-1}H_0''(f)\Bigr).
\end{eqnarray}
The Brieskorn lattice is invariant under $\partial_\tau^{-1}$,
which is well defined as an isomorphism
$\paa_\tau^{-1}:V^{>-1}\to V^{>0}$. Thus $H_0''(f)$
is a free $\C\{\{\paa_\tau^{-1}\}\}$-module of rank $\mu$.
The Fourier-Laplace transform in theorem \ref{5.2} (a)
can be described algebraically by
\begin{eqnarray}\label{6.20}
\paa_\tau^{-1}\mapsto z,\quad \paa_\tau\mapsto z^{-1},\quad
\tau\mapsto -\paa_{z^{-1}}=z^2\paa_z.
\end{eqnarray}
The Fourier-Laplace transform $FL(H_0''(f))$
is a free $\C\{z\}$-module in $V^{>0}_{(z)}$
(the index $(z)$ indicates that here we use the coordinate $z$)
and is invariant under $z^2\paa_z$. Theorem \ref{t5.2} (a)
and \eqref{6.20} give 
\begin{eqnarray}\label{6.21}
G^{(\alpha)}F^p_{St}H^\infty_\lambda
=es(.,\alpha)(z)^{-1}\left(z^{-(m-p)}
\Gr_{V_{(z)}}^{m-p+\alpha}FL(H_0''(f))\right).
\end{eqnarray}
Theorem \ref{t5.2} (b) and theorem \ref{t4.3} (a) say roughly
\begin{eqnarray}\label{6.22}
&&F^{\bullet}_{St}\textup{ has good isotropy properties
w.r.t. }S\\
&\iff&GF^{\bullet}_{St}\textup{ has good isotropy properties
w.r.t. }L^{sym}_\kappa \nonumber\\
&\iff&\Gr_V^\bullet FL(H_0''(f))
\textup{ has good isotropy properties
w.r.t. }P.\nonumber
\end{eqnarray}

In fact, a stronger compatibility of $FL(H_0''(f))$ and $P$ holds
\cite[Theorem 10.28]{He02}\cite[8.1]{He03}:
\begin{eqnarray}\label{6.23}
P:FL(H_0''(f))\times FL(H_0''(f))\to z^{m+1}\C\{z\},
\end{eqnarray}
and the induced symmetric pairing, which is the $z^{m+1}$-coefficient,
\begin{eqnarray}\label{6.24}
FL(H_0''(f)/z\cdot FL(H_0'')\times FL(H_0''(f))/z\cdot FL(H_0''(f))\to\C
\end{eqnarray}
is nondegenerate. 
Here the $\mu$-dimensional space $FL(H_0''(f)/z\cdot FL(H_0'')$ is the
0-fiber of a canonical extension to 0 of a bundle on $\C^*$ 
dual to a bundle of Lefschetz thimbles.

The tuple 
\begin{eqnarray}\label{6.25}
TEZP(f):=(H^\infty_\Z,L^{nor},V^{mod}_{(z)},P,FL(H_0''(f)))
\end{eqnarray}
(the index $(z)$ in $V^{mod}_{(z)}$ 
indicates that the coordinate on $\C$
is here $z$, not $\tau$) is a TERP-structure of weight $m+1$ 
in the sense of \cite[definition 2.12]{He03}. Because of the 
lattice structure $H^\infty_\Z$ (instead of just the real structure
$H^\infty_\R$) we even call it a TEZP-structure.

Finally, consider two singularities $f(x_0,...,x_m)$ and
$g(x_{m+1},...,x_{m+n+1})$ and their 
Thom-Sebastiani sum $f+g$. The canonical isomorphism
\begin{eqnarray}\label{6.26}
(Ml(f),L^{hnor})\otimes (Ml(g),L^{hnor})
\cong (Ml(f+g),L^{hnor})
\end{eqnarray}
induces a canonical isomorphism
\begin{eqnarray}\label{6.27}
&&\bigcup_{z\in\Delta'}H^m(f^{-1}(z),\Z)\otimes
H^n(g^{-1}(z),\Z)\\
&\to& \bigcup_{z\in\Delta'}
H^{m+n+1}((f+g)^{-1}(z),\Z)\nonumber
\end{eqnarray}
which respects the pairings $P$ 
(the $L$ in \eqref{5.2} is here $L^{nor}$).

The following theorem was essentially shown in 
\cite[Lemma (8.7)]{SS85}.
But see the remarks after it for some critic.

\begin{theorem}\label{t6.4}
The TEZP-structures satisfy for $f$ and $g$ as above 
the following Thom-Sebastiani formula:
\begin{eqnarray}\label{6.28}
TEZP(f+g)\cong TEZP(f)\otimes TEZP(g).
\end{eqnarray}
\end{theorem}

{\bf Proof:} The isomorphism for the data
$(H^\infty_\Z,L^{nor})$ is the classical Thom-Sebastiani result
in \eqref{6.3} and \eqref{6.8}.
The isomorphism for $P$ follows from its definition
with $L^{nor}$ in \eqref{5.2}. 
The isomorphism for $V^{mod}_{(z)}$ is trivial.

The isomorphism
\begin{eqnarray}\label{6.29}
FL(H_0''(f+g))\cong FL(H_0''(f))\otimes FL(H_0''(g))
\end{eqnarray}
is not so difficult to see, when one looks at
the (imitations of) oscillating integrals behind
the sections in $FL(H_0''(f))$.
If $\sigma_1\in H_0''(f)\cap \bigoplus_{-1<\alpha<N}C^\alpha_\tau$
for some arbitrarily large $N$ and if 
$\delta_1(z)\in H_n(f^{-1}(z),\C)$, then
\begin{eqnarray}\label{6.30}
FL(\sigma_1)(\delta(z))
&=&\int_0^{z\cdot(+\infty)}e^{-\tau_1/z}\cdot
\sigma_1(\delta_1(\tau_1)) d\tau_1,\\
\textup{where }\sigma_1(\delta_1(\tau_1))&=&
\sum_{-1<\alpha<N}\sum_{k=0}^n a(\sigma_1,\alpha,k)\cdot 
\tau_1^\alpha \cdot (\log\tau_1)^k
\nonumber\\
\textup{ for some }&&a(\sigma_1,\alpha,k)\in\C,
\nonumber
\end{eqnarray}
and analogously for $FL(\sigma_2)(\delta_2(z))$ 
if $\sigma_2\in H_0''(g)\cap \bigoplus_{-1<\alpha<N}C^\alpha_\tau$
and $\delta_2(z)\in H_m(g^{-1}(z),\C)$.
The construction of the topological isomorphism \eqref{4.8} of 
Milnor lattices in \cite[I.2.7]{AGV88} gives
for $\tau\in ]0,z\cdot(+\infty)[$
\begin{eqnarray}\label{6.31}
(\sigma_1\otimes\sigma_2)((\delta_1\otimes\delta_2)(\tau))
&=& \int_0^\tau\sigma_1(\delta_1(\tau_1))
\cdot\sigma_2(\delta_2(\tau-\tau_1))d\tau_1.\hspace*{1cm}
\end{eqnarray}
We obtain (with $\tau=\tau_1+\tau_2$ in the second
equality)
\begin{eqnarray}
&&FL(\sigma_1\otimes\sigma_2)((\delta_1\otimes\delta_2)(z))
\nonumber \\
&=& \int_0^{z\cdot(+\infty)}e^{-\tau/z}\cdot
(\sigma_1\otimes\sigma_2)((\delta_1\otimes\delta_2)(\tau))d\tau
\nonumber \\
&=& FL(\sigma_1)(\delta_1(z))\cdot
FL(\sigma_2)(\delta_2(z)).\label{6.32}
\end{eqnarray}
This proves the isomorphism \eqref{6.29}.
\hfill $\Box$

\begin{corollary}\label{t6.5}
Steenbrink's Hodge filtration satisfies for $f$ and $g$
as above the following Thom-Sebastiani formula:
\begin{eqnarray}
&&G(F^p_{St}) H^\infty_{e^{-2\pi i\alpha}}(f+g)\nonumber\\
&=& \sum_{\beta,\gamma,q,r:(*)}
G(F^q_{St}) H^\infty_{e^{-2\pi i \beta}}(f) \otimes
G(F^r_{St}) H^\infty_{e^{-2\pi i \gamma}}(g)\label{6.33}
\end{eqnarray}
where $0<\alpha\leq 1$ and
\begin{eqnarray*}
(*)&:& 0<\beta,\gamma\leq 1,\ \beta+\gamma=(\alpha\textup{ or }
\alpha+1),\\
&&(m-q+\beta)+(n-r+\gamma)=m+n+1-p+\alpha.
\end{eqnarray*}
\end{corollary}

{\bf Proof:} Apply theorem \ref{t6.4} and 
\eqref{6.21}.\hfill$\Box$

\begin{remarks}\label{t6.6}
(i) The isomorphism \eqref{6.29} for $H_0''$ was essentially proved in 
\cite[(8.7) Lemma]{SS85}.
Though Scherk and Steenbrink did not make the compatibility with
the topological Thom-Sebastiani isomorphism
between the cohomology bundles precise,
and they avoided the use of the Fourier-Laplace transformation.
They obtained a $\paa_\tau^{-1}$-linear isomorphism
$H_0''(f+g)\cong H_0''(f)\otimes H_0''(g)$.

\medskip
(ii) They applied this isomorphism together with \eqref{6.19}
in order to obtain a Thom-Sebastiani formula for
Steenbrink's Hodge filtration $F^\bullet_{St}$
\cite[Theorems (8.2) and (8.11)]{SS85}:
It is \eqref{6.33} without the twists by $G$.
But in the cases with $N\neq 0$, this twist is necessary,
in these cases their formula is not correct.

In their proof, they mixed 
$\paa_\tau^{-1}$-linearity and $\tau$-linearity.
They extracted from the isomorphism 
$H_0''(f+g)\cong H_0''(f)\otimes H_0''(g)$ maps
$C^\beta(f)\otimes C^\gamma(g)\to C^{\beta+\gamma}(f+g)$
in the variable $\tau$ \cite[Lemma (8.8)]{SS85}
and went with them into the defining
formula \eqref{6.19} of $F^{\bullet}_{St}$.

Of course, in the case $N=0$, the isomorphism
$G$ in definition \ref{t4.2} is just a rescaling,
and then $G(F^\bullet_{St})=F^\bullet_{St}$,
so then their Thom-Sebastiani formula is correct.

\medskip
(iii)  
In the case $g=x_{m+1}^2$, the sum $f+g=f+x_{m+1}^2$ is a
suspension of $f$. Theorem \ref{t6.4} gives in that case 
an isomorphism
\begin{eqnarray}\label{6.34}
TEZP(f+x_{n+1}^2)\cong TEZP(f)\otimes TEZP(x_{n+1}^2),
\end{eqnarray}
and formula \eqref{6.33}
boils down to theorem \ref{t4.6}.

\medskip
(iv)
We expect that the following generalization of theorem \ref{t4.6}
holds: {\it For any two (signed or not) Steenbrink PMHS, the  
formula \eqref{6.32} gives a (signed or not) Steenbrink PMHS.}
\end{remarks}

\end{document}